\documentclass[12pt,reqno,NoRunningHeads]{amsart}
\usepackage{latexsym}
\usepackage{amsmath}
\usepackage{amssymb}
\usepackage{amsfonts}
\usepackage{verbatim}
\usepackage[latin1]{inputenc}
\usepackage{color}

\newtheorem{Th}{Theorem}[section]

\newtheorem{cor}[Th]{Corollary}

\newtheorem{proposition}[Th]{Proposition}

\newtheorem{Cor}[Th]{Corollary}

\newtheorem{Def}[Th]{Definition}

\newtheorem{rk}[Th]{Remark}

\let \ssection=\section
\renewcommand{\section}{\setcounter{equation}{0}\ssection}

\def\^#1{\if#1i{\accent"5E\i}\else{\accent"5E #1}\fi}
\def\"#1{\if#1i{\accent"7F\i}\else{\accent"7F #1}\fi}

\newcommand{\RR}{\mathbb R}

\newcommand{\<}{\langle}
\renewcommand{\>}{\rangle}
\def\^#1{\if#1i{\accent"5E\i}\else{\accent"5E #1}\fi}
\def\"#1{\if#1i{\accent"7F\i}\else{\accent"7F #1}\fi}
\allowdisplaybreaks%[2]

\begin{document}
\title[Rate of Convergence of Implicit Approximations]
{Rate of Convergence of Implicit Approximations\\
for stochastic evolution equations}
\author[I. Gy\"ongy]
{Istv\'an Gy\"ongy}
\thanks{This paper was written while the first named author was visiting
the University of Paris 1. The research of this author is partially supported by
EU Network HARP}
\address{School of Mathematics and Maxwell Institute for Mathematical Sciences, 
University of Edinburgh,
King's  Buildings,
Edinburgh, EH9 3JZ, United Kingdom}
\email{gyongy@maths.ed.ac.uk}

 \author[A. Millet]{Annie Millet}
\thanks{The research of the second named author is partially supported
by the research project BMF2003-01345}
\address
{Laboratoire de Probabilit\'es et Mod\`eles Al\'eatoires
(CNRS UMR 7599), Universit\'es Paris~6-Paris~7, Boite Courrier 188, 
4 place Jussieu, 75252 Paris Cedex 05,
{\it  and }\\
\indent Centre d'Economie de la Sorbonne (CNRS UMR 8174), \'Equipe 
 SAMOS-MATISSE,
 Universit\'e Paris 1 Panth\'eon Sorbonne,
 90 Rue de Tolbiac, 75634 Paris Cedex 13}
\email{amil@ccr.jussieu.fr {\it and} annie.millet@univ-paris1.fr} 
%\date{\today}



\begin{abstract}
Stochastic evolution equations 
in Banach spaces 
with unbounded nonlinear drift and diffusion operators are 
considered. Under some regularity condition assumed for the solution,  
the rate of
convergence of implicit Euler approximations 
is estimated under strong monotonicity and
Lipschitz conditions. The results are applied to a class of 
quasilinear stochastic PDEs of parabolic type. 
\end{abstract}

\maketitle
\section{Introduction}
                                            \label{intro}
Let 
$V\hookrightarrow  H\hookrightarrow   V^*$ 
be a
{\it normal triple} of spaces with dense and
continuous embeddings, where $V$
is a separable and reflexive Banach space,
$H$ is a Hilbert space, 
identified with its dual by 
means %the help
of the inner product in $H$, 
and $V^*$ is the dual of $V$.  
Thus $\<v,h\>=(v,h)$ for all 
$v\in V$ and $h\in H^*=H$, 
where $\<v,v^*\>=\<v^*,v\>$ 
denotes the duality product 
of $v\in V$, $v^*\in V^*$, 
and $(h_1,h_2)$ denotes the 
inner product of $h_1,h_2\in H$. Let
$W=\{W(t):t\geq 0\}$ be a 
$d_1$-dimensional 
Brownian motion carried by 
a stochastic basis
$(\Omega, \mathcal F, ({\mathcal F}_t)_{t\geq 0}, P)$. 
Consider the stochastic evolution equation
\begin{equation}                                      \label{u}
 u(t)=u_0 + \int_0^t A(s,u(s))\, ds 
+  \sum_{k=1}^{d_1} \int_0^t
 B_k(s,u(s))\, dW^k(s)\, ,
\end{equation}
where $u_0$ is a $V$-valued 
${\mathcal F}_0$-measurable random variable,
$A$ and $B$ are  (non-linear) 
adapted  operators  defined on
$[0,\infty[\times V\times \Omega$ 
with values in $V^*$ and
$H^{d_1}:=H\times...\times H$, respectively.

It is well-known, see \cite{KR},
\cite{Pa} and
\cite{R},   that this equation admits a 
unique solution if the following conditions are 
met: 
There exist constants $\lambda>0$, $K\geq0$ 
and an ${\mathcal F}_t$-adapted non-negative locally integrable 
stochastic process $f=\{f_t:t\geq0\}$ such that

(i) (Monotonicity) 
There exists a constant $K$ such that  
\begin{equation}                                     \label{monotonicity}
2\<u-v,A(t,u)-A(t,v)\>
+\sum_{k=1}^{d_1}|B_k(t,u)-B_k(t,v)|_H^2
\leq K|u-v|_H^2,
\end{equation}

(ii) (Coercivity) 
$$
2\<v,A(t,v)\>+\sum_{k=1}^{d_1}|B_k(t,v)|_H^2
\leq-\lambda |v|_V^2+K|v|_H^2+f(t), 
$$

(iii) (Linear growth) 
$$
|A(t,v)|_{V^*}^2\leq K|v|_{V^*}^2+f(t),
$$

(iv) (Hemicontinuity) 
$$
\lim_{\lambda\to0}\<w, A(t,v+\lambda u)\>=\<w, A(t,v)\>
$$
hold for all for $u,v, w\in V$, $t\in[0,T]$ and $\omega\in\Omega$. 
 
Under these conditions equation \eqref{u} 
has a unique solution 
$u$ on $[0,T]$.  
(See  Definition \ref{sol} below  for the definition of the
solution.) Moreover,  if
$E|u_0|_H^2<\infty$ and
$E\int_0^Tf(t)\,dt<\infty$, then 
$$
E\sup_{t\leq T}|u(t)|_H^2+
E\int_0^T|u(t)|_V^2\,dt<\infty. 
$$
In \cite{GyMi} it is shown that under these conditions,  
approximations defined by various
implicit and explicit schemes converge to $u$. 

Our aim is to prove rate of convergence 
estimates for these approximations. To achieve this aim 
we require stronger assumptions: a strong monotonicity 
condition on $A,B$ and a Lipschitz condition on $B$ in $v\in V$. In 
the present paper we consider implicit time discretizations. 
Note  that without space discretizations, in general,  explicit
time discretizations do not converge.  Consider, for example, the 
heat equation $du(t)=\Delta u(t)$, with initial condition 
$u(0)=u_0\in L_2(\mathbb R^d)$. Then the explicit time discretization 
on the grid 
$\{k/n\}_{k=0}^n$ gives the approximation $u_n(k/n):=(I+\Delta/n)^ku_0$ 
at time $t=k/n$. Hence clearly, if 
$u_0\notin\cap_{i=1}^{\infty}W^i_2(\mathbb R^d)$, then  
$u(k/n)$ does not belong to the Sobolev space $W^l_2(\mathbb R^d)$, with 
any fixed negative index $l$, when $k$ is sufficiently large.
  
The study of various space-time discretization schemes will be
done in the continuation of the present paper. 

We require also the following  time  
regularity from the solution $u$ (see condition (T2)): 
$E|u_0|_V^2<\infty$, almost surely $u_t\in V$ for every $t\in T$, and 
there exist some constants $C$ and $\nu>0$ such that 
\[ 
 E|u(t)-u(s)|_V^2
\leq C\, |t-s|^{2\nu}\, ,\]
for all $s,t\in[0,T]$. 
Note that unlike the solutions to stochastic differential equations,  
the solutions to stochastic PDEs can satisfy this 
condition with a variety of exponents 
$\nu$, different from $1/2$,  due to the interplay 
between space and time regularities of the solutions. 
(See \cite{Krylov} for space and time regularity of the solutions 
to stochastic parabolic PDEs of second order.)  
Note also that our general setting allows us to cover a large class 
of stochastic parabolic PDEs of order $2m$ for any $m\geq1$ 
(see \cite{KR} for the class of stochastic parabolic SPDEs of 
order $2m$ and see \cite{CW} for the stochastic Cahn-Hilliard equation). 

In the case of 
time independent operators $A$ and $B$ we obtain  
the rate of convergence  
for the implicit approximation 
$u^{\tau}$ corresponding 
to the mesh size $\tau=T/m$ of 
the partition of $[0,T]$
$$
E\max_{i\leq m}|u(i\tau)-u^{\tau}(i\tau)|_H^2
+E\sum_{i\leq  m}|u(i\tau)-u^{\tau}(i\tau)|_V^2\tau
\leq C\tau^{\nu},
$$
 where $C$ is a constant 
independent of $\tau$. If in addition to the above 
assumptions 
$A$ is also Lipschitz continuous in $v\in V$ then the order of 
convergence is doubled, 
$$
E\max_{i\leq  m}|u(i\tau)-u^{\tau}(i\tau)|_H^2
+E\sum_{i\leq  m}|u(i\tau)-u^{\tau}(i\tau)|_V^2\tau
\leq C\tau^{2\nu}.
$$
 In the case of time dependent  
$A$ and $B$ it is natural to assume that they are H\"older continuous 
in $t$ in order to control the error due to their discretization 
in time. However, it is possible to control this discretization error 
when the operator $A$ is not even continuous in $t$, if we 
discretize it by taking the average of $A(s)$ over the intervals 
$[t_i,t_{i+1}]$. This explains the discretization of $A(t)$ 
and condition (T1) below. 
If both operators $A$ and $B$ are 
H\"older continuous in time  then we use also the obvious discretization:      
$A^\tau_{t_i}= A(t_{i+1},.)$ and $B^\tau_{k,t_i}=B(t_i,.)$.  

As examples we present a class of 
quasi-linear  stochastic partial 
differential equations (SPDEs) of parabolic type, 
and show that it  
satisfies our assumptions. Thus 
we obtain rate of convergence results 
also for implicit approximations of linear parabolic 
SPDEs, in particular, for the Zakai 
equation of nonlinear filtering. We refer to 
\cite{KrylovRozovski}, \cite{Par}, \cite{P} and
\cite{R} for basic results for the 
stochastic PDEs of nonlinear filtering.  

We will extend these results to degenerate 
parabolic SPDEs, and to space-time explicit and  
implicit schemes for stochastic evolution equations 
in the continuation of this paper.

In Section 2 we give a precise description 
of the schemes and state the assumptions on the
coefficients which ensure 
the convergence of these schemes to the solution $u$ of (\ref{u}).
In Section 3 estimates for the speed 
of convergence of time implicit schemes are stated and
proved. Finally, in the last section, we give a class of
examples of quasi-linear stochastic PDEs
for which all the assumptions of the main theorem, 
Theorem \ref{convimplunif}, 
are fulfilled.
\medskip

As usual, we denote by $C$ a constant 
which can change from line to line.

\section{Preliminaries and the approximation scheme}

Let $(\Omega, {\mathcal  F}, ({\mathcal F}_t)_{t\geq 0}, P)$ be a
stochastic basis, satisfying the usual conditions, i.e., 
$({\mathcal F}_t)_{t\geq0}$ is an increasing right-continuous 
family of sub-$\sigma$-algebras of $\mathcal F$ such that  
${\mathcal F}_0$ contains every $P$-null set. Let 
$W=\{W(t): t\geq0\}$ be a $d_1$-dimensional  
Wiener martingale with respect to 
$({\mathcal F_t})_{t\geq 0}$, i.e., $W$ is an 
${\mathcal F}_t$-adapted Wiener process 
with values in ${\mathbb R}^{d_1}$ such that  
$W(t)-W(s)$ is independent of ${\mathcal F}_s$ for all 
$0\leq s\leq t$. 

Let $T$ be a given positive number. 
Consider the stochastic evolution equation 
\eqref{u} for $t\in[0,T]$ in a triplet of spaces 
$$
V\hookrightarrow H\equiv H^*\hookrightarrow V^*,
$$
satisfying the following conditions:
$V$ is a separable and reflexive Banach space 
over the real numbers, embedded 
continuously and densely into a Hilbert space 
$H$, which is identified with its dual 
$H^*$ by means 
of the inner product $(\cdot,\cdot)$ in $H$, such that 
$(v,h)=\<v,h\>$ for all $v\in V$ and $h\in H$, where 
$\<\cdot,\cdot\>$ denotes the duality product between $V$ and 
$V^*$, the dual of $V$.  Such triplet of spaces is called 
a normal triplet.

Let us state now our assumptions on the initial value $u_0$ and 
the operators $A$, $B$ in the equation.
Let
\[ A:[0,T] \times V\times \Omega \rightarrow V^*\,
, \quad B:[0,T]\times V\times \Omega \rightarrow H^{d_1}\]
be such that for every $v,w\in V$ and 
$1\leq k\leq d_1$, $\langle A(s,v), w \rangle $ and
$(B_k(s,v),w)$ are adapted processes and
the following conditions hold:
\medskip

{\bf (C1)} The pair $(A,B)$ satisfies the
{\it strong monotonicity condition},   i.e.,
there exist constants  
$\lambda>0$ and $L>0$ such 
almost surely 
$$
2\, \langle u-v, A(t,u)-A(t,v)\rangle
+ \sum_{k=1}^{d_1} |B_k(t,u)-B_k(t,v)|^2_H
$$
\begin{equation}                            \label{monotone}
 + \lambda\, |u-v|_V^2 
\leq L\, |u-v|_H^2  
\end{equation}
for all  $t\in ]0,T]$,
 $u$ and $v$ in $V$.

{\bf (C2)} ({\it Lipschitz  condition on $B$})
There exists a constant $L_1$ such 
that almost surely 
\begin{equation}                                      \label{lipB}
\sum_{k=1}^{d_1} |B_k(t,u)-B_k(t,v)|_H^2 
\leq L_1\, |u-v|_V^2 \,
\end{equation}
for all $t\in [0,T]$, $u$ and $v$ in $V$.

{\bf (C3)} ({\it Lipschitz   condition on $A$})
There exists a constant $L_2$  
 such that
almost surely
\begin{equation}                                            \label{lipA}
|A(t,u)-A(t,v)|_{V^*}^2 \leq L_2\,|u-v|_V^2
\end{equation}
for all $t\in [0,T]$, $u$ and $v$ in $V$.

{\bf (C4)} $u_0:\Omega\rightarrow V$ is
 ${\mathcal F}_0$-measurable and 
 $E|u_0|_V^2 <\infty$. There exist non-negative random variables 
$K_1$ and $K_2$ such that 
$EK_i<\infty$,  
and 
\begin{equation}
\sum_{k=1}^{d_1} |B_k(t,0)|_H^2 \leq K_1           \label{B0}  
\end{equation}
\begin{equation}
|A(t,0)|_{V^*}^2 \leq  K_2                         \label{A0}
\end{equation}
for all $t\in[0,T]$ and $\omega\in\Omega$. 
\smallskip

\begin{rk}                                         \label{19.11.06.06}
If $\lambda=0$ in \eqref{monotone} then 
one says that $(A,B)$ satisfies the 
{\it monotonicity condition}. Notice that 
this condition together with the Lipschitz 
condition \eqref{lipA} on $A$ implies 
the Lipschitz condition \eqref{lipB} on $B$. 
\end{rk}

\begin{rk}
(1) Clearly, \eqref{lipA}--\eqref{A0} 
and \eqref{lipB}--\eqref{B0} imply 
that $A$ and $B$ satisfy the 
{\it  growth condition}  
\begin{equation}                                   \label{growthB}  
\sum_{j=1}^{d_1}|B_k(t,v)|^2_H\leq 2L_1
|v|_{V}^2 + 2K_1,
\end{equation}
 and
\begin{equation}                                   \label{growthA} 
|A(t,v)|_{V^*}^2 \leq 
2L_2 \,|v|_V ^2 +  \,
  2K_2 
\end{equation} 
respectively, for all $t\in[0,T]$, 
$\omega\in\Omega$ and $v\in V$.

(2)
Condition \eqref{lipA} obviously implies that 
the operator $A$ is 
{\it hemicontinuous}:
\begin{equation}                                   \label{hemi}
\lim_{\varepsilon \rightarrow 0}   
\langle A(t,u+\varepsilon v), w\rangle
= \langle A(t,u),w\rangle \, 
\end{equation}
for all  $t\in [0,T]$ and   $u,v,w\in V$.

(3) The strong monotonicity condition (C1), (C2) and
\eqref{B0}, \eqref{A0} yield that the pair $(A,B)$ satisfies  the 
following 
{\it coercivity condition}: there exists a
 non-negative 
random variable  $K_3$     
such $EK_3<\infty$ and almost surely
\begin{equation}        \label{coercive}
2\, \langle v,A(t,v)\rangle  + \sum_{k=1}^{d_1} |B_k(t,v)|^2_H
+ \tfrac{\lambda}{2}\, |v|^2_V \leq L|v|^2_H+  K_3 
\end{equation}
for all $t \in  ]0,T]$, $\omega\in\Omega$ and $v\in V$.
\end{rk}
\begin{proof} We show only (3). By the strong monotonicity 
condition 
\begin{equation}
 2\, \langle v,  A(t,v) \rangle  
+  \sum_{k=1}^{d_1}  |B_k(t,v)|^2_H  
 + \tfrac{\lambda}{2}\, |v|^2_V 
\leq  
L|v|_H^2+R_1(t) +R_2(t) 
\end{equation}
 with
\begin{eqnarray*}
R_1(t)&= & 2\, \langle v,A(t,0)\rangle, \\
R_2(t)&=& \sum_{k=1}^{d_1} |B_k(t,0)|_H^2
 + 2 \sum_{k=1}^{d_1} \Big( B_k(t,v)-B_k(t,0)\, ,\,B_k(t,0) \Big)\, . 
\end{eqnarray*}
Using (C2) and (\ref{A0}), we have  
\begin{align*}
|R_1|&\leq  \tfrac{\lambda}{4}|v|_V^2 + \tfrac{4K_2}{\lambda} \, ,\\
|R_2|&\leq 2\, \left(\sum_{j=1}^{d_1} |
B_k(t,v)-B_k(t,0)|_H^2\right)^{\frac{1}{2}}
\left(\sum_{k=1}^{d_1} |B_k(t,0)|_H^2\right)^{\frac{1}{2}}+ K_1 \\
&\leq \tfrac{\lambda}{4} |v|_V^2 +  C K_1. 
\end{align*}
Thus, (C1) concludes the proof of (\ref{coercive}).

\end{proof}

\begin{Def}                                                             \label{sol}
An $H$-valued adapted continuous process 
$u=\{u(t):t \in[0,T ] \}$ is a solution to
equation \eqref{u} on $[0,T]$ if almost surely  
$u(t)\in V$ for almost every $t\in[0,T]$, 
\begin{equation} \int_0^T  |u(t)|_V^2\, dt <\infty\, , \label{*}\end{equation} 
and 
\begin{equation}\label{solu}
( u(t), v)  =(u_0,v)
+ \int_0^t \langle A(s,u(s)),v
\rangle \, ds
+ \sum_{k=1}^{d_1} \int_0^t (B_k(s,u(s)),v)\,
dW^k(s)
\end{equation}
holds for all $t\in [0,T]$ and $v \in V$. 
We say that the solution to \eqref{u} on $[0,T]$ is unique 
if for any solutions $u$ and $v$ to \eqref{u} on $[0,T]$ 
we have 
$$
P(\sup_{t\in[0,T]}|u(t)-v(t)|_H>0)=0.
$$
\end{Def}

The following theorem is well-known (see 
\cite{KR}, \cite{Pa} and \cite{R}). 

\begin{Th}                                                    \label{existu} 
Let  $A$ and $B$ satisfy 
the monotonicity, coercivity, linear growth 
and hemicontinuity conditions  (i)-(iv)   
formulated in the Introduction. 
Then for every $H$-valued ${\mathcal F}_0$-measurable 
random variable $u_0$,  equation (\ref{u}) has a
unique solution $u$ on $[0,T]$.  
Moreover, if $E|u_0|_H^2<\infty$ 
and $E\int_0^Tf(t)\,dt<\infty$, 
then 
\begin{equation}                                              \label{bounduH}
E\big( \sup_{t\in[0,T]}|u(t)|^2_H\big)
+E\int_0^T  |u(t)|^2_V\,dt <\infty
\end{equation}
holds. 
\end{Th}

Hence by the previous remarks  
we have the following corollary. 

\begin{Cor}                                                       \label{eu}
Assume that conditions  
(C1), (C2) hold.  
Then for every $H$-valued random variable 
$u_0$ equation \eqref{u} 
has a unique solution $u$, and if $E|u_0|_H^2<\infty$, 
then 
\eqref{bounduH} holds.
\end{Cor}

\medskip

\smallskip
\noindent
{\bf Approximation scheme.}
For a fixed integer $m\geq1$ and $\tau:=T/m$ 
we define the approximation $u^\tau$ for the solution $u$  
by an implicit time discretization of equation (\ref{u})
as follows: 
\begin{eqnarray}                                                
u^\tau(t_0)&=& u_0\, ,                                        \nonumber
\\
u^\tau(t_{i+1})&= &u^\tau(t_i)
+ \tau\, 
A^\tau_{t_i}\big(u^\tau(t_{i+1})\big)                          \nonumber
\\
&&+\sum_{k=1}^{d_1}  
{B}^{\tau}_{k,t_i} \big(u^\tau(t_i)\big)\,
\big(W^k(t_{i+1})- W^k({t_i})\big)\, \quad
\text{\rm for $0\leq i<m$},
                                                              \label{time2}  
\end{eqnarray}
where $t_i:=i\tau$ and  
\begin{equation}                                                  \label{Am}
A^\tau_{t_i}(v)= \frac{1}{\tau} \, 
\int_{t_i}^{t_{i+1}}
A(s,v)\, ds \;,
\end{equation}
\begin{equation}                                               \label{tildeBm}       
{B}^{\tau}_{k,0}(v)=0,\,\quad  B^{\tau}_{k,t_{i+1}}(v)
=\frac{1}{\tau}\int_{t_{i}}^{t_{i+1}}B_k (s,v) \, ds                 
\end{equation}
for $i=0,1,2,...,m$.
\smallskip
 
A random vector 
$u^{\tau}:=\{u^{\tau}(t_i):i=0,1,2,...,m\}$ 
is called a solution to 
scheme \eqref{time2} 
if $u^{\tau}(t_i)$ is a $V$-valued 
${\mathcal F}_{t_i}$-measurable 
random variable such that 
$E|u^{\tau}(t_i)|_V^2<\infty$ and \eqref{time2} 
hold for every $i=0, \cdots, m-1$.

We use the notation  
\begin{equation}
                                                             \label{kappa1}
\kappa_1(t):=i\tau
\mbox{ \rm{ for }} \, t\in [i\tau ,(i+1)\tau[, \;
\mbox{\rm and }\;
\kappa_2(t):=(i+1)\tau\,
\mbox{  \rm{  for }} \, t\in ]i\tau,(i+1)\tau]
\end{equation}
for integers $i\geq0$, and set  
$$
A_t(v)=A_{t_i}(v), \quad B_{k,t}(v)=B_{t_i}(v)
$$ 
for $t\in[t_i,t_{i+1}[$, $i=0,1,2,...m-1$ and $v\in V$.

Another possible choice is 
\begin{equation} \label{choice2}
A^\tau_{t_i}(u)=A(t_{i+1},u) \quad \mbox{\rm and} \quad B^\tau_{k,t_i}(u)=B_k(t_i,u)\quad
 \mbox{\rm for}\; i=0, 1, \cdots, m-1.
\end{equation} 

The following theorem establishes the 
existence and uniqueness
of $u^\tau$ 
for large enough $m$, and provides estimates 
in $V$ and in $H$. 
We remark that in practice \eqref{time2} 
should also be solved numerically. This is possible for example 
by Galerkin's approximations and by finite elements methods. 
In the continuation of this paper we consider explicit and implicit 
time discretization schemes together with simultaneous 
`space discretizations',  
and we estimate the error of the corresponding  approximations
 for \eqref{u}. 

\begin{Th}                                      \label{existsnm}
Assume that $A$ and $B$ satisfy the monotonicity, 
coercivity, linear 
growth and hemicontinuity conditions (i)--(iv). 
Assume also that (C4) holds. Let $A^\tau$ and
$B^\tau$ be defined either by \eqref{Am} and \eqref{tildeBm}, 
or by \eqref{choice2}.
Then there exist an integer $m_0$ and a constant 
$C$, such that for  $m\geq m_0$
 equation (\ref{time2}) has a unique solution 
$\{u^\tau(t_i):i=0,1,...,m\}$, and  
\begin{equation}                                       \label{estimHVtau}
E\max_{0\leq i\leq m}\big| u^{\tau}(i\tau)\big|_H^2 
+  E\sum_{i=1}^m\big| u^{\tau}(i\tau)\big|^2_V\tau \,
      \leq C\, .                             
\end{equation}
\end{Th} 

\begin{proof}
For the sake of simplicity, we only give 
the proof in the case $A^\tau$ and $B^\tau$ are
defined by \eqref{Am} and \eqref{tildeBm}.  
This theorem with estimate 
\begin{equation}                                                     
                                                      \label{13.29.05.06}
\max_{0\leq i\leq m}E\big| u^{\tau}(i\tau)\big|_H^2 
+  E\sum_{i=1}^m\big| u^{\tau}(i\tau)\big|^2_V\tau \,
\leq C                                         
\end{equation}
in place of \eqref{estimHVtau} 
is proved in \cite{GyMi} for a slightly different 
implicit scheme. For the above implicit scheme 
the same proof can be repeated without essential 
changes. 
For the convenience of the reader we recall 
from \cite{GyMi} 
that the existence and uniqueness 
of the solution $\{u^{\tau}(t_i):i=0,1, 2,...,m\}$ 
to \eqref{time2}--\eqref{tildeBm} is based on the following 
proposition (Proposition 3.4 from \cite{GyMi}): 
Let  $D:V\rightarrow V^*$ be a mapping  such that

(a) $D$ is monotone, i.e., for every $x,y\in V$,
$\langle D(x) -D(y),x-y\rangle  \geq 0$;

(b) $D$ is hemicontinuous, i.e., 
${\displaystyle \lim_{\varepsilon\rightarrow 0}
  \langle D(x+\varepsilon y),z\rangle  =\langle D(x),z\rangle }$
for every $x,y,z\in V$;

(c) there exist positive constants 
$K$, $C_1$ and $C_2$,  such that 
\begin{equation*}                                    \label{11.04.11}
|D(x)|_{V^*}\leq K\,  (1+|x|_V), 
\quad
\langle D(x),x\rangle \geq C_1\, |x|_V^2- C_2\; , 
\quad \forall x\in V. 
\end{equation*}

Then for every $y\in V^*$, there exists 
$x\in V$ such that $D(x)=y$ and
\begin{equation*}                                        \label{norm}
|x|_V^2 \leq \frac{C_1 +2\, C_2}{C_1} 
+ \frac{1}{C_1^2}\, |y|_{V^*}^2\, .
\end{equation*}

If there exists a positive constant $C_3$ such that
\begin{equation}                                                       \label{strongmono}
 \langle D(x_1)-D(x_2),x_1-x_2\rangle 
\geq C_3\, |x_1-x_2|^2_{V^*}\; ,
\quad \forall x_1,x_2\in V\, ,
\end{equation}
then for any $y\in V^*$, the equation 
$D(x)=y$ has a unique solution $x\in V$.

\noindent
Note that for each $i=1,2,...m-1$ equation 
\eqref{time2} for $x:=u^\tau(t_{i+1})$ 
can be rewritten as $Dx=y$ with 
$$
D:=I-\tau A^\tau_{t_i}, \quad
y:=u^{\tau}(t_{i})+\sum_{k=1}^{d_1}  
{B}^{\tau}_{k,t_i} \big(u^\tau(t_i)\big)\,
\big(W^k(t_{i+1})- W^k({t_i})\big)                          
$$
where $I$ denotes the identity on $V$. It is easy to 
verify that due to conditions (i)--(iv) and 
(C4) the operator $D$ satisfies the conditions  (a), (b) 
and (c)  
for sufficiently 
large $m$. Thus a solution $\{u^{\tau}(t_i):i=0,1,...,m\}$ 
can be obtained by recursion on $i$ for all $m$ greater 
than some $m_0$. To show the uniqueness 
we need only verify  \eqref{strongmono}. By \eqref{Am} 
and by the monotonicity 
condition (i) we have 
$$
\langle D(x_1)-D(x_2),x_1-x_2\rangle=|x_1-x_2|^2_H-
\int_{t_i}^{t_{i+1}}\langle A(s,x_1)-A(s,x_2),x_1-x_2\rangle\,ds
$$
$$
\geq |x_1-x_2|^2_H-K\tau|x_1-x_2|^2_H=(1-K\tau)|x_1-x_2|^2_H,
$$
where the constant $K$ is from \eqref{monotonicity}.   
Hence it is clear that \eqref{strongmono} holds  
if $m$ is sufficiently large.

Now we show \eqref{13.29.05.06}. 
From the definition of $u^{\tau}(t_{i+1})$ we have 
\begin{equation}                                                             \label{16.29.05.06}
|u^{\tau}(t_j)|_H^2=|u_0|^2_H+{\mathcal I}(t_j)+
{\mathcal J}(t_j)+{\mathcal K}(t_j)
-\sum_{i=1}^j|A^{\tau}_{t_i}(u^{\tau}(i\tau))|_H^2\tau
\end{equation}
for $t_j=j\tau$, $j=0,1,2,...m$, where 
\begin{eqnarray*}
{\mathcal I}(t_j)&:=&
2\int_0^{t_j}
\<u^{\tau}(\kappa_2(s)), A(s,u^{\tau}(\kappa_2(s))) \>\,ds, 
\\
{\mathcal J}(t_j)&:=&\sum_{1\leq i<j}\,
|\sum_k 
B^{\tau}_{k,t_i}(u^{\tau}(i\tau)) (W^k(t_{i+1})-W^k(t_i))|_H^2, 
\\
{\mathcal K}(t_j)&:=& 2 \; 
{\sum_k} \int_0^{t_j}
\big( u^{\tau}(\kappa_1(s)),B^{\tau}_{k,s} 
(u^{\tau}(\kappa_1(s))) \big)\,dW^k(s), 
\end{eqnarray*}
and $\kappa_1$, $\kappa_2$ 
are piece-wise constant functions defined 
by \eqref{kappa1}.   
By It\^o's formula  for every $k,l=1,2,...,d_1$
$$
(W^k(t_{i+1})-W^k(t_{i}))
(W^l(t_{i+1})-W^l(t_{i}))
$$
$$
=\delta_{kl}(t_{i+1}-t_i)+
M^{kl}(t_{i+1})-M^{kl}(t_{i}), 
$$
where  
$\delta_{kl}=1$ for $k=l$ and $0$ otherwise, 
and 
$$
M^{kl}(t):=\int_{0}^{t}
\big( W^k(s)-W^k(\kappa_1(s) \big))\,dW^l(s)
+\int_0^t\big( W^l(s)-W^l(\kappa_1(s)) \big)\,dW^k(s). 
$$
Thus we get 
$$
{\mathcal J}(t_j)={\mathcal J}_1(t_j)+{\mathcal J}_2(t_j),
$$
with
$$
{\mathcal J}_1(t_j):=\sum_{1\leq i<j}\sum_k
|{B}^{\tau}_{k,t_i}  (u^{\tau}( t_i))|_H^2\tau
$$
$$
{\mathcal J}_2(t_j):=\int_0^{t_j} 
\sum_{k,l}(B_{k,s}^{\tau}(u^{\tau}(\kappa_1(s))),B_{l,s}^{\tau} 
(u^{\tau}(\kappa_1(s))))\,dM^{kl}(s).  
$$
By the Davis inequality we have 
\begin{align*}
E\max_{j\leq m}|{\mathcal J}_2(t_j)|=
\\
&\leq 3
\sum_{k,l}
E
\left\{
\int_0^T 
|B_{k,s}^{\tau}(u^{\tau}(\kappa_1(s)))|^2_H
|B_{l,s}^{\tau}(u^{\tau}(\kappa_1(s)))|^2_H\,d\<M^{kl}\>(s)
\right\}^{1/2}
\\
&\leq C_1
\sum_{k,l}
E
\left\{
\int_0^T 
|B_{k,s}^{\tau}(u^{\tau}(\kappa_1(s)))|^4_H
|W^l(s)-W^l(\kappa_1(s))\big|^2\, ds
\right\}^{1/2}
\\
& \leq C_1\sum_{k,l}
\, E \Big[ \max_j \big|{B}^{\tau}_{k, t_j}(u^{\tau}(t_j))
\big|_H \; \sqrt{\tau} \\
& \qquad \times  
\; \Big\{ \frac{1}{\tau} \, 
\int_0^T|{B}_{k,s}^{\tau}(u^{\tau}(\kappa_1(s)))|^2_H\,
 \big|W^l(s)-W^{l}(\kappa_1(s))\big|^2 ds \Big\}^{1/2} \Big] 
\\
&  \leq d_1C_1\sum_k \, \tau E  \max_j 
\big|B^{\tau}_{k,t_j}(u^{\tau}(t_j)) \big|_H^2 \\
&\quad + C_1\tau^{-1}\sum_{k,l}\,
E\int_0^T  |{B}_{k,s}^{\tau}(u^{\tau}(\kappa_1(s)))|^2_H\,
 \big|W^l(s)-W^{l}(\kappa_1(s))\big|^2 ds  
\\
&\leq C_2\Big(1+E\sum_{ j \leq 1}\, 
|u^{\tau}(j\tau)|_V^2\tau\Big),
\end{align*}
where $C_1$ and $C_2$ are constants, 
independent of $\tau$. 
Here we use that by
Jensen's inequality for every $k$ 
$$
\sum_{1\leq i<j}
| {B}^{\tau}_{k,t_i}  
(u^{\tau}(i\tau))|_H^2\tau
\leq \int_0^{t_j}|B_k(s,u^\tau(\kappa_2(s))|_H^2 \,ds, 
$$
and that the coercivity 
condition (ii) and the growth condition on (iii)  
imply the growth condition  
\eqref{growthB} on $B$ with some constant $L_1$ and 
random variable $K_1$ satisfying 
$EK_1<\infty$. 
Hence by taking into account the coercivity condition 
we obtain 
\begin{align}
E\; \max_{j\leq m}  \; &  
\big[ {\mathcal I}(t_j)+{\mathcal J}(t_j) 
\big]  
\nonumber \\
&\leq 
E\max_{j\leq m}\int_0^{t_j} \Big[ 2\;  
\big\langle u^\tau (\kappa_2(s))\, ,\, 
A(s, u^{\tau} (\kappa_2(s)))\big\rangle +
\sum_k|B_k(s,u^\tau 
(\kappa_2(s))|_H^2 \Big] \,ds
\nonumber \\
&\quad+E \; \max_{j\leq m}\;
|{\mathcal J}_2(t_j)| 
\nonumber \\
                                                 \label{15.29.05.06}
& \leq C\Big( 
1+ \max_{j\leq m} \;  E|u^{\tau} (j\tau)|_H^2
+E\sum_{j=1}^m|u^{\tau}(j\tau)|_V^2\tau
\Big) 
\end{align}
with a constant $C$ independent of $\tau$. 
By using the Davis 
inequality again we obtain 
\begin{align}
E \; \max_{j\leq m}\; & \big|{\mathcal K}(t_j)\big|
\leq 6\;  E  \left\{\int_0^{T}
\sum_k \big| \big(u^{\tau} (\kappa_1(s)), {B}^{\tau}_{k,s} 
\big( u^{\tau}(\kappa_1(s)) \big)\big|^2\,ds
\right\}^{1/2}
                                                \nonumber 
\\
& \leq 6\; E
\left[ 
\max_{j\leq m} \big|u^{\tau}(j\tau)\big|_H\; 
\left\{ \int_0^{T}
\sum_k
|{B}^{\tau}_{k , s}(u^{\tau}(\kappa_1(s)))|^2_H\,ds
\right\}^{1/2} 
\right] 
                                                 \nonumber 
\\
& \leq \tfrac{1}{2} \; E \max_{j\leq m}| 
u^{\tau} (j\tau)|^2_H
+18\; E\int_0^{T}
\sum_k
\big | {B}^{\tau}_{k,s} 
(u^{\tau}(\kappa_1(s)))\big|^2_H\,ds
\nonumber \\
                                               \label{17.29.05.06}
& \leq \tfrac12\; E \; \max_{j\leq m} \;  |
u^{\tau}(j\tau)|^2_H +C\; \left( 1+E\; 
\sum_{j\leq m } \; |u^{\tau}(j\tau)|^2_V\tau\right)
\end{align}
with a constant $C$ independent of $\tau$. 
From \eqref{13.29.05.06}--\eqref{17.29.05.06} we get 
\begin{align*} 
E \;\max_{j\leq m} & \; | u^{\tau}
(j\tau)|_H^2\leq  E|u_0|^2 +E\; \max_{j\leq m} \; 
\big( {\mathcal I}(t_j)+{\mathcal J}(t_j)
\big) 
+E\;\max_{j\leq m} \; |{\mathcal K}(t_j)|
\\
& \leq \tfrac{1}{2}\; E\; 
\max_{j \leq m} | 
u^{\tau} (j\tau)|^2_H +C\; (1+ \max_{j\leq m} \; 
E|u^{\tau}(j\tau)|_H^2 +E\; 
\sum_{j\leq m} 
|u^{\tau}(j\tau)|^2_V\tau)
\\
& \leq \tfrac{1}{2}\; E \;\max_{j\leq m}\; 
| u^{\tau} (j\tau)|^2_H+C\; (1+L)<\infty
\end{align*}
by virtue of \eqref{13.29.05.06},   which proves the 
estimate \eqref{estimHVtau}.  
\end{proof}

\section{Convergence  results}
In order to obtain a speed of convergence, 
we require further properties from 
$B(t,v)$ and from the solution $u$ of \eqref{u}. 

We assume that there exists a constant $\nu\in]0,1/2]$
such that:

{\bf (T1)} The coefficient $B$ satisfies 
the following {\it time-regularity}:
There exists a constant $C$ 
and a random variable $\eta\geq0$ with 
finite first moment, such that  almost surely
\begin{equation}                                      \label{HolderB} 
\sum_{k=1}^{d_1} 
|B_k(t,v)-B_k(s,v)|_H^2  
\leq|t-s|^{2\nu}(\eta+C|v|_V^2)
\end{equation}
for all $s\in [0,T]$ and $v\in V$.
\smallskip

{\bf (T2)} The solution $u$ to equation 
\eqref{u}  satisfies the following
{\it regularity} property: almost surely $u(t)\in V$ for all $t\in[0,T]$, and  
there exists a constant $C>0$ such that   
\begin{equation}                                   \label{reguimp} 
 E |u(t)-u(s)|^2_V \leq C\, |t-s|^{2\nu}\,
\end{equation}
for all $s,t\in[0,T]$. 

\begin{rk}Clearly, 
\eqref{reguimp} implies 
\begin{equation}                                  \label{10.19.06.06}
\sup_{t\in[0,T]}E|u(t)|^2_V<\infty. 
\end{equation}
\end{rk}

Finally, in order to prove a convergence result in the $H$ norm uniformly in time, we
also have to require the following uniform estimate on the $V$-norm of $u$:

{\bf (T3)} There exists a random variable $\xi$ such that 
$E\xi^2<\infty$ and 
$$
 \sup_{t\leq T}|u(t)|_V   \leq \xi \quad \text{\rm (a.s.)}.
$$
   
\medskip
In order to establish the rate of convergence of the 
approximations we first suppose that the coefficients 
$A$ and $B$ satisfy the Lipschitz
property.
\begin{Th}                                               \label{cvtimeimpl}
Suppose that the conditions (C1)-(C4),    
(T1) and (T2)  hold. Let $A^\tau$ and $B^\tau$
be defined by \eqref{Am} and \eqref{tildeBm}. 
Then there exist a constant $C$ and 
an integer $m_0\geq1 $ such 
that 
\begin{equation}                                          \label{cvtimpl}
\sup_{0\leq l\leq m} E|u(l\tau)-u^\tau(l\tau)|_H^2 +
 E\sum_{j=0}^m|u(j\tau)-u^\tau(j\tau)|_V^2\tau
\leq C\, \tau^{2\nu}\, 
\end{equation}
for all integers $m\geq m_0$. 
\end{Th}
The following proposition plays 
a key role in the proof. 
\begin{proposition}Assume assumptions (i) through (iv) 
from the Introduction  
and let $A^\tau$ and $B^\tau$ be defined by \eqref{Am} and 
\eqref{tildeBm}. 
Suppose, moreover condition (C4).  Then 
\begin{align}
|u(t_l) & -u^\tau(t_l)|_H^2  
=2\, \int_0^{t_l} \big\langle u(\kappa_2(s)) -
u^\tau(\kappa_2(s)) ,
 A(s, u(s)) -  A(s,  u^\tau(\kappa_2(s))) 
\big\rangle\, ds                                        \nonumber 
\\
 &+ \sum_{i=0}^{l-1} 
\left| \int_{t_i}^{t_{i+1}}
\sum_{k=1}^{d_1} \big[B_k(s,  u(s))-{B}_{k,s}^\tau
(u^\tau(t_i)) \big] \, dW^k(s) 
\right|^2_H                                              \nonumber  
\\
&+2\sum_{k=1}^{d_1} \int_0^{t_l} 
\left(B_k(s,  u(s)) -
{B}_{k,s}^\tau (u^\tau(t_i)) \, ,\, u(\kappa_1(s)) -
u^\tau(\kappa_1(s)) 
\right) dW^k(s)                                           \nonumber
\\  
&-  \sum_{i=0}^{l-1} 
\left| 
\int_{t_i}^{t_{i+1}} \big[ A(s,u(s)) - A(s,u^\tau
(t_{i+1}))\, \big]\, ds 
\right|_H^2
                                                          \label{diffimplicite}
\end{align}
\end{proposition}
holds for every $l=1,2,...,m$.
\begin{proof}
Using (\ref{time2})  we have for any $i=0,\cdots, m-1$
\begin{align*}
|u(t_{i+1}) & - u^\tau(t_{i+1})|_H^2  - 
| u(t_i) - u^\tau(t_i)|_H^2 = \nonumber \\
&\quad 2\, \int_{t_i}^{t_{i+1}} 
\big\langle u(t_{i+1}) - u^\tau(t_{i+1}) ,
 A(s, u(s)) -  A(s,  u^\tau(t_{i+1}))
\big\rangle\, ds  \nonumber \\
&\quad + 2\,  \sum_{k=1}^{d_1}  \left( \int_{t_i}^{t_{i+1}} 
\big[ B_k(s,  u(s)) -
{B}_{k,s}^\tau (u^\tau(t_i)) \big] \, dW^k(s) 
\, ,\, u(t_{i+1})-u^\tau(t_{i+1}) \right)                 \nonumber 
\\
&\quad  - \Big| \int_{t_i}^{t_{i+1}} 
\big[ 
A(s,u(s)) - A(s,u^\tau (t_{i+1}))\, 
\big]\, ds                                               \nonumber 
\\
&\qquad\qquad  +
\sum_{k=1}^{d_1}  \int_{t_i}^{t_{i+1}} 
\big[  B_k(s,  u(s)) -
{B}_{k,s}^\tau (u^\tau(t_i)) 
\big] \, dW^k(s) 
\Big|^2_H                                                 \nonumber  
\\
= & 2\, \int_{t_i}^{t_{i+1}} 
\big\langle 
u(t_{i+1}) - u^\tau(t_{i+1}) ,
 A(s, u(s)) -  A(s,  u^\tau(t_{i+1}))
\big\rangle\, ds                                           \nonumber 
\\
&\quad + \left| \sum_{k=1}^{d_1}  
\int_{t_i}^{t_{i+1}} 
\big[  B_k(s,  u(s)) -
{B}_{k,s}^\tau (u^\tau(t_i))
\big] \, dW^k(s) 
\right|^2_H                                         \nonumber 
\\ 
&\quad + 2\, \sum_{k=1}^{d_1} 
\left(\int_{t_i}^{t_{i+1}}  
\big[  B_k(s,  u(s)) -
{B}_{k,s}^\tau (u^\tau(t_i)) \big] \, dW^k(s)
\, , \, u(t_i)-u^\tau(t_i)
\right)   
\\
&\quad -\left| \int_{t_i}^{t_{i+1}} 
\big[ A(s,u(s)) - A(s,u^\tau (t_{i+1}))\, 
\big]\, ds 
\right|_H^2
\end{align*}
Summing up for $i=1,\, \cdots,\,  l-1$, we obtain 
\eqref{diffimplicite}. 
\end{proof}

\medskip
\noindent
{\it Proof of Theorem \ref{cvtimeimpl}.}

\noindent
Taking expectations in both sided of \eqref{diffimplicite}
 and using the strong monotonicity condition (C1), 
we deduce that for $l=1,\, \cdots,\, m$,
$$
E|u(t_l) -u^\tau(t_l)|^2_H                           
$$
\begin{align}
&\leq E \int_0^{t_l} 2  
\big\langle u(\kappa_2(s)) - u^\tau(\kappa_2(s)) ,
A (s,u(\kappa_2(s))) - A (s, 
u^\tau(\kappa_2(s)))\big\rangle ds 
                                                      \nonumber 
\\ 
& \qquad  +
\sum_{k=1}^{d_1} E \int_0^{t_{l-1}} 
|B_k(s,u(\kappa_2(s)))  
- B_k (s,u^\tau(\kappa_2(s)))|_H^2\, ds
+\sum_{k=1}^3 R_k                                    \nonumber \\
 & \leq - \lambda  \, E\int_0^{t_l}  
|u(\kappa_2(s))-u^{\tau}(\kappa_2(s))|_V^2\,
ds + L \tau E |u(t_l)-u^{\tau}(t_l)|_H^2              \nonumber \\  
&\qquad  
+ L \tau \,  \sum_{i=1}^{l-1} 
E  |u(t_i)-u^{\tau}(t_i)|_H^2\, ds +\sum_{k=1}^3
R_k\, , 
\label{majoimptemps}
\end{align}
where
\begin{align*}
R_1=& E\int_0^{t_l} 2\, 
\big\langle u(\kappa_2(s)) - u^\tau(\kappa_2(s)) ,
A(s,u(s)) -A(s,u(\kappa_2(s)))
\big\rangle\, ds\, ,\\ R_2=&
\sum_{k=1}^{d_1} E \int_0^{\tau} |B_k (s,u(s))|_H^2\, ds\, ,
\\
R_3=&\sum_{k=1}^{d_1} 
\sum_{i=1}^{l-1} E \Big[\int_{t_i}^{t_{i+1}} ds \,  
\Big|B_k(s,u(s)) - \frac{1}{\tau}\int_{t_{i-1}}^{t_i} 
B_k (t,u^\tau(t_i) )  \, dt \Big|_H^2              
\\ &\qquad\qquad  
- \int_{t_{i-1}}^{t_i} |B_k (t,u(t_i)) -
B_k (t,u^\tau(t_i))|_H^2\,  dt \Big] .
\end{align*}
The Lipschitz property of $A$ imposed in (\ref{lipA}), 
(\ref{reguimp}) and Schwarz's inequality imply 
\begin{align}
&| R_1|\leq L_2 \, E\int_0^{t_l} | u(\kappa_2(s))
- u^\tau(\kappa_2(s))|_V \,
|u(s) - u(\kappa_2(s))|_V\, ds\,,                    \nonumber
\\ 
&\quad \leq L_2  \left(
E\int_0^{t_l}
\!\!  | u(\kappa_2(s)) - u^\tau(\kappa_2(s))|_V^2
 ds\right)^{\frac{1}{2}}
\left( 
E\int_0^{t_l}\!\!  | u(s) - u(\kappa_2(s))|_V^2
ds
\right)^{\frac{1}{2}}                                \nonumber
\\ 
&\leq \frac{\lambda}{3}
E\int_0^{t_l}  
|u(\kappa_2(s)) - u^\tau(\kappa_2(s))|_V^2 ds + C
\tau^{2\nu}. 
\label{majoR1timeimp}
\end{align}
A similar computation based on (\ref{lipB}) yields
\begin{align*}
|R_3|&\leq\sum_{k=1}^{d_1}\sum_{i=1}^{l-1} 
E \int_{t_{i-1}}^{t_i} dt \frac{1}{\tau}  
\int_{t_i}^{t_{i+1}} ds
\Big( |B_k(s,u(s)) -B_k(t,u^\tau(t_i))|_H^2 
\\ 
&\qquad\qquad  -   |B_k (t,u(t_i))) -
B_k (t,u^\tau(t_i)))|_H^2 \Big)                      \nonumber 
\\ 
&\leq
\frac{\lambda}{3} E\int_0^{t_{l-1}}  
| u(\kappa_2(t)) - u^\tau(\kappa_2(t))|_V^2 dt 
+ C\, R'_3 
\,                                                  \label{majoR3timeimp1}
\end{align*}
where 
$$R'_3=\sum_{k=1}^{d_1}E \frac{1}{\tau}
\int_{t_1}^{t_l} ds 
\int_{\kappa_1(s)-\tau}^{\kappa_1(s)} dt \,  
 |B_k (s,u(s)) -B_k (t, u(\kappa_2(t)))|_H^2.
$$ 
Hence, using
(\ref{lipB}), (\ref{HolderB}) and (\ref{reguimp})
we have 
\begin{align*}
R'_3&\leq\sum_{k=1}^{d_1} E\frac{1}{\tau} \int_{t_1}^{t_l}
ds \int_{\kappa_1(s)-\tau}^{\kappa_1(s)} dt
\Big[ |B_k (s,u(s))- B_k(t,u(s))|_H^2 
\\
&\qquad\quad + |B_k (t,u(s)) -B_k (t,
u(t))|_H^2
 + |B_k (t,u(t))- B_k (t,u(\kappa_2(t)))|_H^2\Big]
\\ 
&\leq
E\int_{t_1}^{t_l} \tau^{2\nu} \, |u(s)|_V^2 \, ds 
+ C E\frac{1}{\tau} 
\int_{0}^{t_{l-1}} dt \int_{\kappa_2(t)}^{\kappa_2(t)+\tau} ds   
\Big[ 
|u(s)-u(t)|_V^2 
\\
\qquad\quad & \qquad \qquad 
+ |u(t)-u(\kappa_2(t)|_V^2\Big] \leq C\, \tau^{2\nu}. 
\end{align*}
Hence 
\begin{equation}                                         \label{majoR3timeimpl2}
|R_3|\leq C\, \tau^{2\nu}+  
\frac{\lambda}{3}
E\int_0^{t_l} |u(\kappa_2(s))-u^\tau(\kappa_2(s))|_V^2 ds\, .
\end{equation}
Furthermore \eqref{growthB}  and \eqref{10.19.06.06} imply 
\begin{equation}                                        \label{11.19.06.06}
|R_2|\leq C\tau
\end{equation}
with a constant $C$ independent of $\tau$. 
 By inequalities    
(\ref{majoimptemps})--(\ref{11.19.06.06}), 
for sufficiently large $m$,
\begin{align}
E|u(t_l) - u^\tau(t_l)|_H^2 &+ 
\frac{\lambda}{3}
E \int_0^{t_l}
|u(\kappa_2(s))-u^\tau(\kappa_2(s))|_V^2 ds
                                                                \nonumber
\\ 
& \leq \sum_{i=1}^{l-1} L\, \tau\, 
E|u(t_i)-u^\tau(t_i)|_H^2 + C\tau^{2\nu}.                \label{gronwaltimeimp}
\end{align}
Since $\sup_m\sum_{i=1}^m L\, \tau  <+\infty$, 
a discrete version of Gronwall's lemma yields
that there exists $C>0$ such that for $m$ large enough 
$$
\sup_{0\leq l\leq m} E|u(t_l)-u^\tau(t_l)|_H^2
\leq C\tau^{2\nu}.
$$
This in turn with (\ref{reguimp}) 
implies  
$$
E\int_0^T |u(s)-u^\tau(\kappa_2(s))|_V^2 ds
\leq C\tau^{2\nu}, 
$$ 
which completes the proof of the theorem.\hfill $\Box$

\medskip

Assume now that the solution $u$ of equation \eqref{u} satisfies 
also the condition (T3) 
Then we can improve the estimate \eqref{cvtimpl} 
in the previous theorem. 

\begin{Th}                                                  \label{convimplunif} 
Let (C1)-(C4) and (T1)--(T3) hold, 
and let $A^\tau$ and $B^\tau$ be defined by
\eqref{Am} and \eqref{tildeBm}.  
Then for all sufficiently large $m$ 
\begin{equation}                                            \label{19.29.05.06}
E\max_{0\leq j\leq m}|u(j\tau)-u^\tau(j\tau)|_H^2 +
 E \sum_{j=0}^m |u(j\tau)-u^\tau(j\tau)|_V^2\tau
\leq C\, \tau^{2\nu} 
\end{equation}
holds, where $C$ is a constant independent of $\tau$. 
\end{Th}

\begin{proof}   
For $k=1, \cdots, d_1$, set 
\[ F_k(t)=B_k(t,u(t)) - {B}^\tau_{k,t}(u^\tau(\kappa_1(t)))
\]
$$
m(t)=\sum_{k=1}^{d_1}\int_{0}^t F_k(s)\, dW^k(s)\quad 
\mbox{\rm  and }\quad G(s)=m(s)-m(\kappa_1(s))
$$
Then by It\^o's formula 
$$
|m(t_{i+1})-m(t_{i})|_H^2 
= 2 \int_{t_i}^{t_{i+1}}
  \sum_k(G(s)\, ,\,  F_k(s))\,dW^k(s)
$$
$$
+ \sum_{k=1}^{d_1} \int_{t_i}^{t_{i+1}} |F_k(s)|_H^2\, ds
$$ for $i=0,...,m-1$.
Hence by 
using (\ref{diffimplicite}) 
we deduce that for $l=1, \cdots, m$
\begin{equation}                             \label{24.11.06.06}
|u(t_l)-u^\tau(t_l)|_H^2 \leq I_1(t_l)+I_2(t_l)
+2M_1(t_l)+2M_2(t_l)
\end{equation}
with
$$ 
 I_1(t):=2\int_0^{t} 
\langle 
u(\kappa_2(s)) - u^\tau(\kappa_2(s))\, , 
\, A(s,u(s)) - A(s,u^\tau(\kappa_2(s))) 
\rangle
\, ds, 
$$
$$
I_2(t):=\sum_{k=1}^{d_1} \int_0^{t} 
|B_k(s, u(s)) - {B}_{k,s}^{\tau}(u^\tau(\kappa_1(s)))|_H^2\, ds,
$$
$$
M_1(t):=\sum_{k=1}^{d_1} 
\int_0^{t} 
\big( G(s)\, , \, F_k(s)\big) \, dW^k(s),
$$
$$
M_2(t):= \sum_{k=1}^{d_1} \int_0^{t} 
\big( 
F_k(s)\, , \, u(\kappa_1(s)) - u^\tau(\kappa_1(s)) 
\big)\, dW^k(s). 
$$
By (C3) 
$$
\sup_{0\leq l\leq m} |I_1(t_l)| 
\leq\int_0^T|u(\kappa_2(s))-u^{\tau}(\kappa_2(s))|^2_V\,ds
+L_2\int_0^T|u(s)-u^{\tau}(\kappa_2(s))|^2_V\,ds
$$
$$
\leq (1+2L_2)\sum_{i=1}^m
|u(t_i)-u^{\tau}(t_i)|^2_V\tau
+2L_2\int_0^T|u(s)-u(\kappa_2(s))|^2_V\,ds.
$$ 
Hence by Theorem \ref{cvtimeimpl} 
and by condition (T2)
\begin{equation}                                         \label{ineqsupH1}
E\sup_{0\leq l\leq m} |I_1(t_l)|\leq C\tau^{2\nu},
\end{equation}
where $C$ is a constant independent of $\tau$. 

Using
Jensen's inequality, \eqref{growthB} 
and condition (T1) we
have for $s\leq \tau$
\begin{equation}                                      \label{n1}
 \sum_k |F_k(s)|_H^2=\sum_k |B_k(s,u(s))|_H^2 \leq 2\, L_1\, |u(s)|_H^2 +2\, K_1,
\end{equation}
while for $s\in [t_i,t_{i+1}]$, $1\leq i\leq m$, one has for some constant $C$ independent of $\tau$
\begin{align}                \label{n2}
\sum_k |F_k(s)|_H^2 & \leq \frac{1}{\tau} \sum_k \int_{t_{i-1}}^{t_i}
|B_k(s, u(s)) - B_k(r,u^{\tau}(t_i))|_H^2\, dr
\nonumber \\
&\leq 3 \frac{1}{\tau} \sum_k \int_{t_{i-1}}^{t_i} \Big[ |B_k(s, u(s)) - B_k(r, u(s))|_H^2
\nonumber \\
&\qquad +  |B_k(r, u(s)) - B_k(r, u(t_i))|_H^2 +  |B_k(r, u(t_i))-B_k(u^\tau(t_i))|_H^2\Big]\, dr
\nonumber\\
& \leq C\,\Big[  \tau^{2\nu}\, \big(\eta +|u(s)|_V^2\big) +   |u(s)-u(t_i)|_V^2
+  |u(t_i)-u^\tau(t_i)|_V^2\Big] .
\end{align}
Thus, \eqref{n1} and \eqref{n2} yield
\begin{align*}
\sup_{0\leq l\leq m} |I_2(t_l)| & \leq \;  C\int_0^{\tau} |u(s)|_V^2\, ds + C\tau +
 C\, \tau^{2\nu}\, \int_0^T \big( \eta + C\, |u(s)|_V^2 \big)\, ds \\
& \qquad  + C \, \int_0^T |u(s)-u(\kappa_2(s))|_V^2\, ds
+ C \, \sum_{i=1}^m |u(t_i)-u^\tau(t_i)|_V^2\, \tau.
\end{align*}
Hence by Theorem \ref{cvtimeimpl} 
and by condition (T2)
\begin{equation}                                         \label{ineqsupH2}
E\sup_{0\leq l\leq m} |I_2(t_l)|\leq C\tau^{2\nu},
\end{equation}
where $C$ is a constant independent of $\tau$. 
Since $\sup_{0\leq s\leq T}\sum_k|F_k(s)|^2_H$
need not be measurable, we denote by $\Gamma$ 
the set of random variables $\zeta$  
satisfying 
$$
\sup_{0\leq s\leq T}\sum_k|F_k(s)|^2_H\leq \zeta 
\quad \text{\rm (a.s.)}.
$$
For $\zeta \in \Gamma$,  the Davis inequality, and the simple inequality 
$ab\leq \tfrac{\tau}{2}a^2+\tfrac{1}{2\tau}b^2$ yield  
\begin{align}
E \sup_{1\leq l\leq m} &|M_1(t_l)| 
\leq  3\, E 
\left( \int_0^T\sum_{k=1}^{d_1} 
| \big( F_k(s)\, , \,  G(s) \big) |^2 \,ds 
\right)^{\frac{1}{2}}                                               \nonumber  
\\ 
&\leq 3\, E\left(\zeta^{1/2}\,
\, 
\left[\int_0^T  |G(s)|_H^2\, ds
\right]^{\frac{1}{2}}\right)                                       \nonumber 
\\
&\leq\frac{3}{2} \tau \inf_{\zeta\in\Gamma}E\zeta  
+  \frac{3}{2\tau}E \int_0^T  |G(s)|_H^2\, ds .
                                                               \label{21.11.06.06}
\end{align}
By (\ref{growthB})  and \eqref{n2} we deduce 
\begin{align*}
\sup_{0\leq s\leq T} \sum_k|F_k(s)|_H^2 &  
\leq
C\, \tau^{2\nu}\, \big( \sup_{0\leq s\leq T} 
|u(s)|_V^2 + \eta + 1 \big)
 +C\, \max_{1\leq i\leq m} |u(t_i)-u^\tau(t_i)|_V^2
\\
& \leq C \,\Big( 1+\xi +  \max_{1\leq i\leq m} |u(t_i)-u^\tau(t_i)|_V^2\Big) ,
\end{align*}
where $\xi$ is the random variable from
 condition (T3) and $C$ is a constant, independent
 of $\tau$. Hence Theorem  \ref{cvtimeimpl} yield
\begin{equation}                    \label{16.06.06.1}
  \tau \inf_{\zeta \in \Gamma} E\zeta \leq \tau\, C\, (E\eta + E\xi) +
 C\, \tau\, \sum_{i=1}^m E|u(t_i)-u^\tau(t_i)|_V^2
\leq C_1\, \tau^{2\nu},
\end{equation}
where  $C_1$ is a constant, independent
of $\tau$.
Similarly, due to conditions (T1)-(T2)
and Theorem \ref{cvtimeimpl}
\begin{align} \label{14.06.01}
E\sum_k\int_0^T|F_k(s)|_H^2\,ds  & \leq
C\, \tau^{2\nu} \Big( 1+E\int_0^T |u(s)|_V^2\, ds \Big) \nonumber  \\
&\qquad + C\, \tau^{2\nu} +C\, \tau \,
E\sum_{i=1}^m |u(t_i)-u^\tau(t_i)|_V^2 \leq C\, \tau^{2\nu}
\end{align}
with a constant $C$, independent of $\tau$.
Furthermore, the isometry of stochastic
integrals  and (\ref{ineqsupH4}) yield
\begin{align}
\frac{1}{\tau} E\int_0^T |G(t)|_H^2 \, dt 
& \leq  \frac{1}{\tau} E\int_0^T
\left| 
\int_{\kappa_1(t)}^t\sum_k F_k(s)\, dW^k(s) 
\right|_H^2\, dt \nonumber 
\\
& \leq \frac{1}{\tau} E\int_0^T dt 
\int_{\kappa_1(t)}^t \sum_k|F_k(s)|_H^2\, ds 
\leq C\, \tau^{2\nu}.
                                                         \label{ineqsupH5}
\end{align}
Thus from \eqref{21.11.06.06} 
by \eqref{16.06.06.1} and \eqref{ineqsupH5} we have 
\begin{equation}                                           \label{ineqsupH3}
E \sup_{1\leq l\leq m} |M_1(t_l)|\leq C\tau^{2\nu}
\end{equation}
Finally, the Davis inequality implies
\begin{align}
E\sup_{1\leq l\leq m}& |M_2(t_l)|_H\leq 3 \,
E\left( \int_0^T 
\sum_k|\big( F_k(s)\, ,\, u(\kappa_1(s))-u^\tau(\kappa_1(s)) \big)
|^2\, ds
\right)^{\frac{1}{2}}                                        \nonumber  
\\
&\leq \frac{1}{4} \, E\sup_{1\leq l\leq m} 
|u(\kappa_1(s))-u^\tau(\kappa_1(s)) \big)|_H^2 +
18 \, E\int_0^{t_j} |F_k(s)|_H^2\, ds.
                                                            \label{ineqsupH4}
\end{align}
Thus, from 
\eqref{24.11.06.06} by inequalities 
(\ref{ineqsupH1}), (\ref{ineqsupH2}),
(\ref{ineqsupH3}) and (\ref{ineqsupH4}) 
we obtain 
\[ \frac{1}{2} \, E\sup_{1\leq l\leq m} 
|u(t_l) -u^\tau(t_l)|_H^2 \leq C\, \tau^{2\nu},\]
with a constant $C$, independent of $\tau$, 
which with (\ref{cvtimpl}) 
completes the proof of the theorem.
\end{proof}

We now prove that if the coefficient $A$  
does  not satisfy the Lipschitz property 
(C3) but only the
coercivity and growth conditions 
(\ref{growthA})-(\ref{coercive}), then the
order of convergence is divided by two.

\begin{Th}                                                \label{convtimeimp2}
Let $A$ and $B$ satisfy the conditions 
(C1), (C2) and (C4).  Suppose
that conditions (T1) and (T2) hold, 
and let $A^\tau$ and $B^\tau$ be defined by
\eqref{Am} and \eqref{tildeBm}. 
Then there exists a constant $C$, independent of 
$\tau$, such that for all sufficiently large $m$ 
\begin{equation}                                          \label{vitconvtimimp2}
\sup_{0\leq j \leq m}E|u(j\tau )-u^\tau(j\tau)|_H^2 
+ E \sum_{j=1}^m 
|u(j\tau)-u^\tau(j\tau)|_V^2 \, \tau  \leq C\, \tau^{\nu}\, .
\end{equation} 
\end{Th}

\begin{proof}
Using  (\ref{diffimplicite}), 
taking expectations  and using (C1) with $u(s)$ 
and $u^\tau(\kappa_2(s))$,
 we obtain for every $l=1\, \cdots, m$
\begin{align}                                             \label{majoBistimeimp}
E|u(t_l)-u^\tau(t_l)|_H^2 & 
\leq -\lambda E \int_0^{t_l}  
|u(s)-u^\tau(\kappa_2(s)|_V^2\, ds \nonumber\\
&\quad + E \int_0^{t_l} K_1 \, 
|u(s)-u^\tau(\kappa_2(s)|_H^2\, ds 
+ \sum_{k=1}^3 \bar{R}_i\, ,
\end{align}
where
\begin{align*}
\bar{R}_1=&\sum_{j=1}^r 
2E\int_0^{t_l} \big\langle u(\kappa_2(s))- u(s)\, ,\,
A(s, u(s)) - A(s, u^\tau(\kappa_2(s)))
 \big\rangle \, ds\, ,
\\
\bar{R}_2=& \sum_{k=1}^{d_1} E\int_0^{\tau} 
|B_k (s,u(s))|_H^2\, ds \, ,\\
\bar{R}_3=& \sum_{k=1}^{d_1}  \sum_{i=1}^{l-1} 
E\frac{1}{\tau} \int_{t_i}^{t_{i+1}} ds \int_{t_{i-1}}^{t_i}\, dt 
\Big[  
|B_k(s,u(s)) - B_k (t,u^\tau(t_i)))|^2_H 
\\ &\qquad \qquad - 
|B_k (t,u(t)) - B_k (t, u^\tau(t_i))|^2_H 
\Big]\, .
\end{align*}
Using (\ref{growthA}), (\ref{reguimp}), (\ref{10.19.06.06}) and 
Schwarz's inequality, we deduce 
\begin{align}
|\bar{R_1}|&\leq C\,  E \int_0^{t_l}  
|u(\kappa_2(s)) - u(s)|_V 
\big[   
|u(s)|_V + |u^\tau(\kappa_2(s))|_V +  K_2  
\big]\, ds
\nonumber \\
&\leq C\, \left( E\int_0^{t_l} \!\!  |u(s)-u(\kappa_2(s))|_V^2\, ds 
 \right)^{\frac{1}{2}} \; 
\left(  E\int_0^{t_l}\!\!  \Big( |u(s)|_V^2 + |u(\kappa_2(s))|_V^2\Big)
 \, ds\right)^{\frac{1}{2}}
\nonumber  \\
&\quad+  C\, 
\left(E\int_0^{t_l}  
|u(s)-u(\kappa_2(s))|_V^2\, ds\right)^{\frac{1}{2}} 
\nonumber \\
& \leq  C\tau^{\nu}\, .\label{majoR1Bistimeimp}
\end{align}
Furthermore, Schwarz's inequality, (C2) 
and computations similar to that proving 
 (\ref{majoR3timeimpl2}) 
yield for any $\delta>0$ small enough
\begin{align*}
|\bar{R}_3|&\leq \delta \;  \sum_{k=1}^{d_1}  
 E\int_0^{t_{l-1}}\!\! | B_k (t,u(t)) 
- B_k (t,u^\tau(\kappa_2(t)))|_H^2 \,  dt
\\
& \quad + C \; \sum_{k=1}^{d_1}  \sum_{i=1}^{l-1}
\frac{1}{\tau}
 E\int_{t_i}^{t_{i+1}}ds \, \int_{t_{i-1}}^{t_i} dt \, 
|B_k (s,u(s))) - B_k (t,u(t))|_H^2 
\\
&\leq \frac{\lambda}{2} E\int_0^{t_{l-1}}\!\!  
|u(s)-u^\tau(\kappa_2(s))|_V^2 \, ds +
C \tau^{2\nu} \, .
\end{align*}
This inequality 
and (\ref{majoR1Bistimeimp})   
imply that 
\begin{align*}
E|u(t_l)-u^{\tau}(t_l)|_H^2 
&+ \frac{\lambda}{2}\, E \int_0^{t_l} 
 |u(s)-u^\tau(\kappa_2(s))|_V^2
ds\\
& \leq K_1 \int_0^{t_l}  
E|u(s)-u^\tau(\kappa_2(s) )|_H^2 \, ds  + C\, \tau^\nu.
\end{align*}
Hence for any $t\in [0,T]$, 
\begin{align*}
E|u(t)&-u^\tau(\kappa_2(t)|_H^2 
\leq 2\, E|u(\kappa_2(t))-u^\tau (\kappa_2(t))|_H^2 + 2 \, E |u(t)
-u(\kappa_2(t))|_H^2\\
&\leq 2\, K_1 \int_0^{\kappa_2(t)} 
E |u(s)-u^\tau(\kappa_2(s))|_H^2 ds + C\, \tau^\nu + 
 2 \, E |u(t) -u(\kappa_2(t))|_H^2\\
&\leq 2\,  K_1\int_0^t 
E|u(s)-u^\tau(\kappa_2(s))|_H^2\, ds + C\, \tau^\nu +  2 \, E |u(t)
-u(\kappa_2(t))|_H^2\\ &\qquad\qquad + C\, \tau\, 
\Big[ \sup_s E\big(
|u(s)|_H^2 + u^\tau(\kappa_2(s))|_H^2\big) 
\Big].
\end{align*}
It\^o's formula and (\ref{coercive}) 
imply that for any $t\in [0,T]$,
\begin{align*}
 E|u(t)&-u(\kappa_2(t))|_H^2 = E\int_t^{\kappa_2(t)} 
\Big[ 2\langle 
A(s, u(s))\,,\, u(s)\rangle +
\sum_{k=1}^{d_1} |B_k(s,u(s))|_H^2
\Big]\, ds\\
&\leq K_1\, E\int_t^{\kappa_2(t)} |u(s)|_H^2\, ds 
\leq K_1\, \tau\, \sup_{0\leq s \leq T} E|u(s)|_H^2.
\end{align*}
Hence (\ref{bounduH}) and (\ref{estimHVtau}) 
imply that
\[ E|u(t)-u^\tau(\kappa_2(t))|_H^2 
\leq 2\, K_1 \int_0^t E|u(s)-u^\tau(\kappa_2(s))|_H^2\, ds  
+ C\, \tau^\nu\]
and Gronwall's lemma yields 
\begin{equation}                                     \label{egy}
\sup_{0\leq t\leq T}
 E|u(t)-u^\tau(\kappa_2(t))|_H^2 \leq C\tau^\nu .
\end{equation} 
Therefore, 
\begin{equation}                                    \label{ketto}
E\int_0^T   |u(t)-u^\tau(\kappa_2(t)|_V^2\, dt 
<C\tau^{\nu}
\end{equation}
follows by \eqref{majoBistimeimp}. Finally taking 
into account that by (T2) there exists a constant $C$ 
such that 
$$
E|u(t)-u(\kappa_2(t)|_V^2
\leq C\tau^{2\nu} \text{\rm for all $t\in[0,T]$}, 
$$
from \eqref{egy} and \eqref{ketto} we obtain 
\eqref{vitconvtimimp2}. 
\end{proof}

Using the above result one can easily obtain 
the following theorem in the same way as 
Theorem \ref{cvtimeimpl} is
obtained from Theorem \ref{convimplunif}.   

\begin{Th}                                             \label{20.29.05.06}
Let $A$ and $B$ satisfy the conditions 
(C1), (C2) and (C4).  Suppose
that conditions (T1)--(T3)  hold and let $A^\tau$ and $B^\tau$
be defined by \eqref{Am} and \eqref{tildeBm}.    
Then there exists a constant $C$ 
such that for $m$ 
large enough,
\begin{equation} \label{vitconvtimimp2bis}
E\max_{0\leq j\leq m}  |u(j\tau )-u^\tau(j\tau)|_H^2 
+ E \sum_{j=0}^m  
|u(j\tau)-u^\tau(j\tau)|_V^2\tau \leq C\, \tau^{\nu}\, .
\end{equation} 
\end{Th}

\begin{rk}
By analyzing their proof, it is not difficult 
to see that Theorems \ref{cvtimeimpl}, 
\ref{convimplunif}, \ref{convtimeimp2} 
and  \ref{20.29.05.06} 
remain true, if instead of 
\eqref{Am} 
and \eqref{tildeBm}, one uses  
\eqref{choice2} in the definition of the implicit scheme, and requires
furthermore that
$A$ satisfies the following time-regularity similar to (T1):
there exist a constant $C\geq0$ and a random variable $\eta\geq0$ 
with finite expectation,  
such that  almost surely 
\[ |A(t,u)-A(s,u)|_{V^*}^2 \leq |t-s|^{2\nu} 
\big( \eta + C \|u\|_V^2\big)
\] 
for $0\leq s\leq t\leq T$ and $u\in V$.
\end{rk}

\section{Examples}

\subsection{Quasilinear stochastic PDEs}
\noindent Let us consider the stochastic 
partial differential equation 
\begin{align}
du(t,x)=& \big( Lu(t)
+F(t,x,\nabla u(t,x),u(t,x)\big)\,dt                    \nonumber 
\\
                                                    \label{20.29.05.06.eq}
 &+\sum_{k=1}^{d_1} 
\big( M_ku(t,x)+G_k(t,x,u(t,x))\big)\,dW^k(t), 
\end{align}
for $t\in(0,T]$, $x\in{\mathbb R}^d$ 
with initial condition 
\begin{equation}                                      \label{21.29.05.06}
u(0,x)=u_0(x), \quad x\in{\mathbb R}^d,
\end{equation}
where $W$ is a $d_1$-dimensional Wiener martingale 
with respect to the filtration $({\mathcal F}_{t})_{t\geq0}$, 
$F$ and $G_k$ are Borel
functions of 
$(\omega,t,x,p,r)\in \Omega\times[0,\infty)
\times{\mathbb R}^d\times{\mathbb R}^d
\times{\mathbb R}$  and of $(\omega,t,x,r)\in
\Omega\times[0,\infty)\times{\mathbb R}^d\times{\mathbb R}$, 
respectively, and 
$L$, 
$M_k$ are differential operators of the  form 
\begin{equation}                                      \label{linop}
L(t)v(x)=\sum_{|\alpha|\leq 1, |\beta|\leq1}
D^{\alpha}(a^{\alpha\beta}(t,x)D^{\beta}v(x)), 
\quad M_k(t)v(x)=\sum_{|\alpha|\leq1}
b^{\alpha}_k(t,x)D^{\alpha}v(x), 
\end{equation}
with functions 
$a^{\alpha\beta}$ and $b^{\alpha}_k$ of 
$(\omega,t,x)\in \Omega\times[0,\infty)\times{\mathbb R}^d$, 
for all multi-indices 
$\alpha=(\alpha_1,...,\alpha_d)$, 
$\beta=(\beta_1,...,\beta_d)$ of length 
$|\alpha|=\sum_i\alpha_i\leq 1$, $|\beta|\leq1$. 

Here, and later on 
$D^{\alpha}$ denotes $D^{\alpha_1}_1...D^{\alpha_d}_d$ 
for any multi-indices 
$\alpha=(\alpha_1,...,\alpha_d)\in\{0,1,2,...\}^d$, 
where 
$
D_i=\frac{\partial}{\partial x_i} 
$
and $D_i^0$ is the identity operator.  

We use the notation   
$\nabla_p:=(\partial/\partial p_1,...,\partial/\partial p_d)$.   
For $r\geq 0$ let $W^r_2({\mathbb R}^d)$ denote the space 
of Borel functions $\varphi:{\mathbb R}^d\to \mathbb R$ 
whose derivatives up to order $r$ are square integrable 
functions. The norm $|\varphi|_r$ of $\varphi$ in
$W^r_2$ is defined by 
$$
|\varphi|_r^2=\sum_{|\gamma|\leq r}
\int_{{\mathbb R}^d}|D^{\gamma}\varphi(x)|^2\,dx.
$$
In particular, 
$W^2_0({\mathbb R}^d)=L_2({\mathbb R}^d$) 
and $ |\varphi|_0:=|\varphi|_{L_2({\mathbb R}^d)}$.  
Let us use the notation 
${\mathcal P}$ for  the $\sigma$-algebra of predictable 
subsets of $\Omega\times[0,\infty)$, and 
${\mathcal B}(\mathcal{\mathbb R}^d)$ for  
the Borel $\sigma$-algebra 
on ${\mathbb R}^d$.

We fix an integer
$l\geq0$ and assume that the following conditions hold. 
\bigskip

\noindent
{\bf Assumption (A1)} 
(Stochastic parabolicity).
{\it 
There exists a constant 
$\lambda>0$ such that 
\begin{equation}                                   \label{22.29.05.06}
\sum_{|\alpha| = 1,|\beta| = 1}
\left( a^{\alpha\beta}(t,x)-\tfrac{1}{2} 
\sum_{k=1}^{d_1}
b^{\alpha}_kb^{\beta}_k(t,x)\right)
\,z^{\alpha}\, z^{\beta}\geq 
\lambda
\sum_{|\alpha|=1}|z^{\alpha}|^2
\end{equation}
for all $\omega\in\Omega$, $t\in[0,T]$, 
$x\in{\mathbb R}^d$ and 
$z=(z^1,...,z^d)\in{\mathbb R}^d$, where 
$z^{\alpha}
:=z_1^{\alpha_1}z_2^{\alpha_2}...z_d^{\alpha_d}$ 
for $z\in{\mathbb R}^d$ and multi-indices 
$\alpha=(\alpha_1,\alpha_2,...,\alpha_d)$. 
} 

\smallskip

\noindent
{\bf Assumption (A2)} (Smoothness of the linear term).
{\it  The derivatives of $a^{\alpha\beta}$ 
and $b^{\alpha}_k$ up to
order
$l$  are  
${\mathcal P}\otimes{\mathcal B}(\mathcal{\mathbb R}^d)$
-measurable real functions 
such that for a constant $K$  
\begin{equation}                                    \label{23.29.05.06}
|D^{\gamma}a^{\alpha\beta}(t,x)|\leq K, 
\quad |D^{\gamma} b^{\alpha}_k(t,x)|\leq K, 
\quad 
\text{\rm for all 
$|\alpha|\leq1$, $|\beta|\leq1$, 
$k=1,\cdots, d_1$,}
\end{equation}                                       
for all $\omega\in\Omega$, $t\in[0,T]$, 
$x\in{\mathbb R}^d$ and multi-indices 
$\gamma$ with $|\gamma|\leq l$. 
}

\smallskip

\noindent
{\bf Assumption (A3)} 
(Smoothness of the initial condition). 
{\it Let $u_0$ be a $W^l_2$-valued 
${\mathcal F}_0$-measurable random variable  
such that 
\begin{equation}                                   \label{25.29.05.06}
E|u_0|^2_l<\infty.
\end{equation}
}

\smallskip
  
\noindent
{\bf Assumption (A4)} (Smoothness of the nonlinear term).
{\it  The function  $F$ and their 
first order partial derivatives in $p$ 
and $r$ are 
${\mathcal P}\otimes{\mathcal B}(\mathcal{\mathbb R}^d)
\otimes{\mathcal B}(\mathcal{\mathbb R}^d)
\otimes{\mathcal B}(\mathcal{\mathbb R})$-measurable 
functions, and 
$g_k$  and  its  
first order derivatives in 
$r$ are 
${\mathcal P}
\otimes{\mathcal B}(\mathcal{\mathbb R}^d)
\otimes{\mathcal B}(\mathcal{\mathbb R})$
-measurable functions for every $k=1,..,d_1$. 
There exists 
a constant $K$ such that 
\begin{equation}                                      \label{24.29.05.06}
|\nabla_p F(t,x,p,r)|
+|\tfrac{\partial}{\partial r} F(t,x,p,r)|
+\sum_{k=1}^{d_1}
|\tfrac{\partial}{\partial r} G_k(r, x)|\leq K
\end{equation}
for all $\omega\in\Omega$, $t\in[0,T]$, 
$x\in{\mathbb R}^d$, $p\in{\mathbb R}^d$ 
and $r\in\mathbb R$.     
There exists a random variable $\xi$ 
with finite first  
moment, such that 
\begin{equation}                                       \label{20.17.06.06}
 |F(t,\cdot,0,0)|_0^2
+\sum_{k=1}^{d_1}|G_k(t,\cdot,0)|_0^2\leq\xi
\end{equation}
for all $\omega\in\Omega$ and $t\in[0,T]$.
}  
\smallskip

\begin{Def}
An $L_2({\mathbb R}^d)$-valued continuous 
${\mathcal F}_t$-adapted process $u=\{u(t):t\in[0,T]\}$ 
is called a generalized
solution  to the Cauchy problem 
\eqref{20.29.05.06.eq}-\eqref{21.29.05.06} 
on $[0,T]$ 
if almost surely $u(t)\in W^1_2(\RR^d)$ for almost every $t$, 
$$
\int_0^T|u(t)|^2_{1}\,dt<\infty, 
$$
and 
\begin{align*}
d(u(t),\varphi)=&
\Big\{
\sum_{|\alpha|\leq1,|\beta|\leq1} 
(-1)^{|\alpha|}\big(a^{\alpha\beta}\, D^{\beta}u(t)\, ,\,
D^{\alpha}\varphi
\big) +\big( F(t,\nabla u(t),u(t))\, ,\, \varphi\big)
\Big\}\,dt 
\\
& + \sum_{k=1}^{d_1} \Big\{ \sum_{|\alpha|\leq 1}
\big(b^{\alpha}_kD^{\alpha}u(t)\, ,\,
\varphi\big)+ 
\big(G_k(t,u(t))\, ,\, \varphi\big) \Big\}\,dW^k(t)
\end{align*}
holds on $[0,T]$ for every 
$\varphi\in C^{\infty}_0({\mathbb R}^d)$, where 
$(v,\varphi)$ denotes the inner product of $v$ 
and $\varphi$ in 
$L_2({\mathbb R}^d)$.  
\end{Def}

Set $H=L_2({\mathbb R}^d)$,  
$V=W^1_2({\mathbb R}^d)$ and consider the normal triplet 
$
V\hookrightarrow H\hookrightarrow V^*
$
based on the inner product in $L_2({\mathbb R}^d)$, 
which determines the duality $\<\,,\,\>$ between $V$ 
and $V^*=W^{-1}_2({\mathbb R}^d)$. 
By \eqref{23.29.05.06}, \eqref{24.29.05.06} and 
\eqref{20.17.06.06} there exist a constant $C$ 
and a random variable $\xi$ with finite first moment, 
 such that 
$$
\Big|\sum_{|\alpha|\leq1,|\beta|\leq1} 
(-1)^{|\alpha|}\big(a^{\alpha\beta}(t)\, D^{\beta}v\, ,\,
D^{\alpha}\varphi
\big)\Big|\leq C|v|_1|\varphi|_1, 
\;\;
\sum_{k=1}^{d_1}
|(b^{\alpha}_k(t) D^{\alpha}v\, ,\,
\varphi\big)|^2
\leq C |v|_0^2|\varphi|_0^2, 
$$
$$
|\big( F(t,\nabla v,v)\, ,\, \varphi\big)|^2  
\leq C|v|_1^2|\varphi|_1^2+\xi,
\quad
\sum_{k=1}^{d_1}|(G_k(t,u(t))\, ,\, \varphi\big)|^2
\leq C|v|_1^2|\varphi|_0^2 +\xi
$$
for all  $\omega$, 
$t\in[0,T]$ and $v,\varphi\in V$. 
Therefore the operators $A(t)$, $B_k(t)$ defined by 
\begin{align}
\<A(t,v),\varphi\>=&
\sum_{|\alpha|\leq1,|\beta|\leq1} 
(-1)^{|\alpha|}\big(a^{\alpha\beta}(t) \, D^{\beta}v\, ,\,
D^{\alpha}\varphi
\big) +\big( F(t,\nabla v,v)\, ,\, \varphi\big),
\nonumber\\
                                            \label{AB}
(B_k(t,v)\,,\,\varphi)=&\; 
\big(b^{\alpha}_k(t) D^{\alpha}v\, ,\,
\varphi\big)+ 
\big(G_k(t,v)\, ,\, \varphi\big), 
\quad v,\varphi\in V  
\end{align}
are mappings from $V$ into $V^*$ and $H$, respectively, 
for each $k$ and $\omega,t$, such that the growth 
conditions \eqref{growthB} and \eqref{growthA} hold. 
Thus we can cast the 
Cauchy problem \eqref{20.29.05.06.eq}--
\eqref{21.29.05.06}  into the evolution equation 
\eqref{u}, and it is 
an easy exercise to show that Assumptions 
(A1), (A2) with $l=0$ and Assumption (A4) 
ensure that conditions 
(C1) and (C2) hold. 
Hence Corollary \ref{eu} gives the following result.  

\begin{Th}                                                \label{01.30.05.06}
Let Assumptions  
(A1)-(A4) hold with $l=0$. 
Then problem \eqref{20.29.05.06.eq}-\eqref{21.29.05.06} 
admits a unique
generalized  solution $u$ on $[0,T]$. Moreover, 
\begin{equation}                                            \label{02.30.05.06}
E\Big( \sup_{t\in[0,T]}|u(t)|_0^2\Big) 
+E\int_0^T|u(t)|^2_1\,dt<\infty.
\end{equation}
\end{Th}

Next we formulate a result on the regularity 
of the generalized solution. We  
need the following assumptions.

\smallskip

\noindent
{\bf Assumption (A5)} 
{\it 
The first order derivatives of $G_k$ in $x$ 
are 
${\mathcal P}
\otimes{\mathcal B}(\mathcal{\mathbb R}^d)
\otimes{\mathcal B}(\mathcal{\mathbb R})$
-measurable functions, and there exist 
a constant $L$, a 
${\mathcal P}\otimes{\mathcal B}(\mathcal{\mathbb R})$
-measurable function $K$ of $(\omega,t,x)$ and a random variable 
$\xi$ with finite first moment, 
such that 
$$
\sum_{k=1}^{d_1}|D^{\alpha}G_k(t,x,r)|
\leq L|r|+K(t,x), 
\quad |K(t)|_0^2\leq \xi  
$$
for all multi-indices $\alpha$ with $|\alpha|=1$, for all 
$\omega\in\Omega$, $t\in[0,T]$, 
 $x\in {\mathbb R}^d$ and $r\in\mathbb R$. 
}

\smallskip

\noindent
{\bf Assumption (A6)} 
{\it 
The first order derivatives of $F$ in $x$ are 
${\mathcal P}\otimes{\mathcal B}(\mathcal{\mathbb R}^d)
\otimes{\mathcal B}(\mathcal{\mathbb R}^d)
\otimes{\mathcal B}(\mathcal{\mathbb R})$-measurable 
functions,  
and there exist a constant $L$, a  
${\mathcal P}\otimes{\mathcal B}(\mathcal{\mathbb R})$
-measurable function $K$ of $(\omega,t,x)$ 
and a random variable 
$\xi$ with finite first moment, 
such that 
$$
|\nabla_xF(t,x,p,r)|\leq L(|p|+|r|)+K(t,x), 
\quad 
|K(t)|_0^2\leq\xi
$$
for all $\omega,t,x, p,r$. 
}

\smallskip

\noindent
{\bf Assumption (A7)} 
{\it 
There exist  
${\mathcal P}\otimes{\mathcal B}(\mathcal{\mathbb R})$
-measurable functions $g_k$ such that 
$$
G_k(t,x,r)=g_k(t,x)
\quad \text{\rm for all $k=1,2,...,d_1$, $t,x,r$}, 
$$
and the derivatives in $x$ of $g_k$ up to order 
$l$ are ${\mathcal P}\otimes{\mathcal B}(\mathcal{\mathbb R})$
-measurable functions such that 
$$
\sum_{k=1}^{d_1}|g_k(t)|_l^2\leq \xi,  
$$
for all $(\omega,t)$, where $\xi$ 
is a random variable with finite first 
moment.
} 

\begin{Th}                                      \label{17.13.06.06}
Let Assume  
(A1)-(A4) with $l=1$. 
Then for the generalized  solution 
$u$ of 
\eqref{20.29.05.06.eq}-\eqref{21.29.05.06}
the following statements hold:

(i) Suppose (A5). Then $u$  
is a $W^1_2({\mathbb R}^d)$-valued 
continuous process and  
\begin{equation}                                   \label{03.30.05.06}
E\Big(\sup_{t\leq T}|u(t)|_1^2\Big)
+E\int_0^T|u(t)|^2_2\,dt<\infty\,;
\end{equation}

(ii) 
Suppose (A6) and (A7) with $l=2$.  
Then $u$ is a $W^2_2({\mathbb R}^d)$-valued 
continuous process and  
\begin{equation}                               \label{07.30.05.06}
E\Big(\sup_{t\leq T}|u(t)|_2^2\Big)
+E\int_0^T|u(t)|^2_3\,dt<\infty\,. 
\end{equation}
\end{Th}

\begin{proof}Define 
$$
\psi(t,x)=F(t,x,\nabla u(t,x),u(t,x)), 
\quad \phi_k(t,x)=G_k(t,x,u(t,x)) 
$$
for $t\in[0,T]$, $\omega\in \Omega$ 
and ${x\in\mathbb R}^d$, 
where $u$ is the generalized solution of
\eqref{20.29.05.06.eq}-\eqref{21.29.05.06}. 
Then due to \eqref{02.30.05.06} 
$$
E\int_0^T|\psi(t)|_0^2\,dt<\infty, 
\quad E\sum_k\int_0^T|\phi_k(t)|_1^2\,dt<\infty. 
$$
Therefore, the Cauchy problem 
\begin{eqnarray}                                 
dv(t,x)&=&(Lv(t,x)+\psi(t,x))\,dt
\nonumber \\
                                                 \label{05.30.05.06}
& & +\sum_{k=1}^{d_1}(M_kv(t,x)+\phi_k(t,x))\,dW^k(t), 
\quad t\in(0,T], 
\quad x\in{\mathbb R}^d\, , \\
                                                  \label{06.30.05.06}
v(0,x)&=& u_0(x), \quad x\in{\mathbb R}^d
\end{eqnarray}
has a unique generalized solution $v$ 
on $[0,T]$. Moreover, 
by Theorem 1.1 from \cite{KrRo},  
$v$ is a $W^1_2$-valued continuous 
${\mathcal F}_t$-adapted 
process and 
$$
E\Big(\sup_{t\leq T}|v(t)|_1^2\Big)
+E\int_0^T|v(t)|^2_2\,dt<\infty. 
$$
Since $u$ is a generalized solution to 
\eqref{05.30.05.06}--\eqref{06.30.05.06}, 
by virtue of the uniqueness of the
generalized  solution we have $u=v$, 
which proves (i).  
Assume now (A6) and (A7). Then obviously  
(A5) holds, and therefore 
due to \eqref{03.30.05.06} 
$$
E\int_0^T|\psi(t)|_1^2\,dt<\infty, 
\quad  
E\sum_k\int_0^T|\phi_k(t)|_2^2\,dt<\infty. 
$$ 
Thus by Theorem 1.1 of \cite{KrRo} 
the generalized solution $v=u$ of 
\eqref{05.30.05.06}--\eqref{06.30.05.06}
is a $W^2_2({\mathbb R}^d)$-valued 
continuous process such that 
\eqref{07.30.05.06} holds. 
The proof of the theorem is complete.  
\end{proof}

\begin{Cor}                                  \label{coregularity} 
Let (A1)-(A4) hold with $l=2$. 
Assume also (A6) and (A7).  
Then there exists a constant $C$ such that for the generalized 
solution $u$ of \eqref{20.29.05.06.eq}--\eqref{21.29.05.06}
we have 
$$
E|u(t)-u(s)|_1^2\leq C|t-s|
\quad \text{\rm for all $s,t\in[0,T]$}.
$$ 
\end{Cor}

\begin{proof} By the theorem on It\^o's formula from 
\cite{KR} (or see \cite{GK:82}) 
from almost surely 
$$
u(t)=u_0+\int_0^t(Lu(s)+\psi(s))\,ds+\sum_{k=1}^{d_1}
\int_0^t(M_ku(s)+g_k(s)\,dW^k(s)
$$
holds, as an equality in $L_2({\mathbb R}^d)$,  
for all $t\in[0,T]$, where  
$$
\psi(s,\cdot) :=F(s,\cdot\nabla u(s,\cdot),u(s,\cdot)). 
$$
Due to (ii) from Theorem \ref{17.13.06.06}
\begin{align*}
E\Big| \int_s^t & \big( Lu(r)+\psi(r) \big)\,dr\Big|_1^2
\leq E\Big( \int_s^t|Lu(r)+\psi(r)|_1\,dr \Big)^2 
\\
& \leq |t-s|\, E\int_s^t |Lu(r)+\psi(r)|_1^2\,dr
\\
& \leq C\, |t-s| \, 
\left( E\int_0^T|u(t)|_3^2\,dt
+E\int_0^T|\psi(t)|_1^2\,dt\right)
\leq C|t-s|
\end{align*}
for all $s,t\in[0,T]$, where $C$ is a constant. Furthermore, 
by Doob's inequality 
\begin{align*}
E\left|\int_s^tM_ku(r)+g_k(r)\,dW^k(r)\right|^2_1 
& \leq 
4\int_s^tE|M_ku(r)+g_k(r)|_1^2\,dr
\\
& \leq C_1|t-s|
\Big[
1+E\Big( \sup_{t\leq T}|u(t)|_2^2\Big)
\Big]\leq C_2|t-s|
\end{align*}
for all $s,t\in[0,T]$, 
where $C_1$ and $C_2$ are constants. 
Hence 
\begin{align*}
E|u(t)-u(s)|_1^2\leq &  2E\,\Big| 
\int_s^t \big( Lu(r)+\psi(r) \big)\,dr
\Big|_1^2
\\ & \quad
+ 2E\Big| \sum_{k=1}^{d_1} \int_s^t
(M_ku(r)+g_k(r))dW^k(r)\Big|_1^2\leq C|t-s|, 
\end{align*}
and the proof of the corollary is complete. 
\end{proof}

\medskip 

The implicit scheme \eqref{time2} applied to problem 
\eqref{20.29.05.06.eq}-\eqref{21.29.05.06} 
reads as follows.  
\begin{eqnarray}                                          
u^\tau(t_0)&=& u_0\,
,                                                           \nonumber 
\\ u^\tau(t_{i+1})&= &u^\tau(t_i)
+ 
\left(L^\tau_{t_i}u^\tau(t_{i+1})
+F^{\tau}_{t_i}(u^{\tau}(t_{i+1})\right)\tau
                                                            \nonumber 
\\
&&+\sum_{k=1}^{d_1}
\left(M^{\tau}_{k,t_i}u^\tau(t_i)
+G^{\tau}_{k,t_i}(u^\tau(t_i))\right)\,
(W^k(t_{i+1})- W^k({t_i}))\, ,                             \label{Ptime2}
\end{eqnarray}                                                             
for $0\leq i<m\,$,   where  
\begin{align}
L^\tau_{t_i}v:& =
\sum_{|\alpha|\leq1,|\beta|\leq1}
D^{\alpha}(a^{\alpha\beta}_{t_i}(x)D^{\beta}v), 
\quad  
M^{\tau}_{k,t_i}:=\sum_{|\alpha|\leq1}
b^{\alpha}_{k,t_i}D^{\alpha}v, 
\nonumber \\
                                                  \label{PAm}
a^{\alpha\beta}_{t_i}(x):& = \frac{1}{\tau} \, 
\int_{t_i}^{t_{i+1}}
a^{\alpha\beta}(s,x)\, ds ,
\\
                                                 \label{PBm}      
{b}^{\alpha}_{k,0}(x) & =0,\,\quad 
b^{\alpha}_{k,t_{i+1}}  (x)
=\frac{1}{\tau}\int_{t_{i}}^{t_{i+1}}b_k (s,x) \, ds,                 
\\
\nonumber F^{\tau}_{t_i}(x,p,r):& =\frac{1}{\tau}
\int_{t_i}^{t_{i+1}}F(s,x,p,r)\,ds, 
\\
G^{\tau}_{k,0}(x,r):& =0, 
\quad
G^{\tau}_{k,t_{i+1}}(x,r):=
\int_{t_i}^{t_{i+1}}G_k(s,x,r)\,ds. \nonumber
\end{align} 
\smallskip

\smallskip

\begin{Def} A random vector 
$\{u^{\tau}(t_i):i=0,1,2,...,m\}$ is a called a 
generalized solution 
of the scheme \eqref{Ptime2} if $u^{\tau}(t_0)=u_0$, 
$u^{\tau}(t_i)$ is a $W^1_2({\mathbb R}^d)$-valued 
${\mathcal F}_{t_i}$-measurable random variable 
such that 
$$
E|u^{\tau}(t_i)|^2_1<\infty
$$
and almost surely 
\begin{align*}
(u^{\tau}(t_i),\varphi)= & 
\sum_{|\alpha|\leq1,|\beta|\leq1}(-1)^{|\alpha|}
(a^{\alpha\beta}_{t_i}D^{\beta}u^{\tau}(t_i),
D^{\alpha}\varphi)\tau
+(F_{t_{i-1}}^{\tau}
(\nabla u_{t_{i-1}}^{\tau},u_{t_{i-1}}^{\tau}),
\varphi)\tau \\
&\;  +\sum_{k}(\sum_{|\alpha|\leq 1}
b^{\alpha}_{t_{i-1}}D^{\alpha}u^{\tau}_{t_{i-1}}
+G_{k,t_{i-1}}(u^{\tau}_{t_{i-1}}),\varphi)
(W^{k}(t_{i})-W^k(t_{i-1}))
\end{align*}
for $i=1,2,...,m$ and all  
$\varphi\in C_0^{\infty}({\mathbb R}^d)$, where 
$(\cdot,\cdot)$ is the inner product in 
$L_2({\mathbb R}^d)$.  
\end{Def} 
From this definition it is clear that, using the operators $A$, $B_k$ 
defined by \eqref{AB},  
we can cast the scheme \eqref{Ptime2} 
into the abstract scheme \eqref{time2}. Thus by applying Theorem 
\ref{existsnm} we 
get the following theorem.  
\begin{Th}                                          \label{kozeli}
Let 
(A1)-(A4) hold with $l=0$. Then there exists an integer 
$m_0$ such that \eqref{Ptime2} has a
unique  generalized solution 
$\{u^{\tau}(t_i):i=0,1,...,m\}$ for every $m\geq m_0$. 
Moreover, there exists a constant $C$ such that 
$$
E\max_{0\leq i\leq m}|u^{\tau}(t_i)|^2_0
+E\sum_{i=1}^m|u^{\tau}(t_i)|^2_1\leq C
$$
for all integers $m\geq m_0$. 
\end{Th}
To ensure condition (T1) to hold we impose 
the following assumption.   

\smallskip

\noindent
{\bf Assumption (H)} {\it There exists a constant $C$ and 
a random variable $\xi$ with finite first moment  
such that for $k=1,2,...,d_1$ 
\begin{align*}
|D^{\gamma}(b^{\alpha}_k(t,x)-b^{\alpha}_k(s,x))|
& \leq C|t-s|^{1/2}\quad
\text{\rm for all}\;  \omega\in\Omega, 
x\in{\mathbb R}^d \;  \text{\rm and } \; 
|\gamma|\leq l, 
\\
 |g_k(s)-g_k(s)|_l^2 & \leq \xi|t-s|
\end{align*}
for all $s,t\in[0,T]$. 
}

\smallskip
Now applying Theorem \ref{convimplunif}
we obtain the following result. 
\begin{Th}                                    \label{fopelda}
Let (A1)-A(4) and (A6)-(A7) hold with $l=2$. 
Assume (H) with $l=0$.  Then 
\eqref{20.29.05.06.eq}--\eqref{21.29.05.06} 
and \eqref{Ptime2} have a unique 
generalized solution $u$ and 
$u^{\tau}=\{u^{\tau}(t_i):i=0,1,2,...,m\}$, 
respectively, for all integers $m$ larger 
than some integer $m_0$. Moreover, for all integers 
$m>m_0$ 
\begin{equation}                               \label{sebesseg}
E\max_{0\leq i
\leq m}|u(i\tau)-u^{\tau}(i\tau)|_0^2
+E\sum_{i=1}^{m}|u(i\tau)-u^{\tau}(i\tau)|^2_1\tau
\leq C\tau,
\end{equation}
where $C$ is a constant, independent of $\tau$.  
\end{Th}
\begin{proof}
By Theorems \ref{01.30.05.06}  and \ref{kozeli} 
the problems 
\eqref{20.29.05.06.eq}--\eqref{21.29.05.06}  and 
\eqref{Ptime2} have 
a unique solution $u$ and $u^\tau$, respectively. 
It is an easy exercise to verify 
that Assumption (H) ensures that condition (T1) 
holds. By virtue of Corollary \ref{coregularity} 
condition (T2) is valid with $\nu=1/2$. 
Condition (T3) clearly holds 
by statement (i) 
of Theorem  \ref{17.13.06.06}. 
Now we can apply Theorem \ref{convimplunif} , 
which gives \eqref{sebesseg}. 
 
\end{proof}
\subsection{Linear stochastic PDEs} 
Let Assumptions (A1)-(A3) and (A7) hold and  
impose also the following condition  
on $F$. 

\smallskip
\noindent
{\bf Assumption (A8)} 
{\it  There exist  a 
${\mathcal P}\otimes{\mathcal B}(\mathcal{\mathbb R})$
-measurable function $f$ such that 
$$
F(t,x,p,r)=f(t,x), \quad \text{\rm for all $t,x,p,r$}, 
$$
and the derivatives in $x$ of $f$ up to order 
$l$ are ${\mathcal P}\otimes{\mathcal B}(\mathcal{\mathbb R})$
-measurable functions such that 
$$
|f(t)|_l^2\leq \xi,  
$$
for all $(\omega,t)$, where $\xi$ 
is a random variable with finite first 
moment.
} 

Now equation \eqref{05.30.05.06} has become the 
linear stochastic PDE  
\begin{equation}                               \label{11.12.06.06}
du(t,x)=
(Lu(t,x)+
f(t,x))\,dt +
\sum_{k=1}^{d_1}(M_ku(t,x)+g_k(t,x))\,dW^k(t),                       
\end{equation}
and by Theorem \ref{convimplunif}  
we have the following result. 

\begin{Th} \label{th4.6} Let $r\geq0$ be an integer. 
Let Assumptions (A1)--(A3) and (A7)--(A8) 
hold with $l:=r+2$, and let Assumption (H) hold 
with $l=r$.   
 Then there is an integer $m_0$ 
such that \eqref{11.12.06.06}--\eqref{21.29.05.06} 
and \eqref{Ptime2} have a unique 
generalized solution $u$ and 
$u^{\tau}=\{u^{\tau}(t_i):i=0,1,2,...,m\}$, 
respectively, for all integers $m>m_0$. 
Moreover, 
\begin{equation}                               \label{lsebesseg}
E\max_{0\leq i
\leq m}|u(i\tau)-u^{\tau}(i\tau)|_r^2
+E\sum_{i=1}^{m}
|u(i\tau)-u^{\tau}(i\tau)|^2_{r+1}\tau
\leq C\tau
\end{equation}
holds for all $m>m_0$, where $C$ is a constant independent 
of $\tau$. 
\end{Th}
\begin{proof}For $r=0$ the statement of this theorem 
follows immediately from Theorem \ref{fopelda}. 
For $r>0$ 
set $H=W_2^{r}({\mathbb R}^d)$ and   
$V=W^{r+1}_2({\mathbb R}^d)$ and consider 
the normal triplet  
$
V\hookrightarrow H\equiv H^*\hookrightarrow V^*
$
based on the inner product 
$(\cdot\,,\cdot):=(\cdot\,,\cdot)_r$ in $W^r_2({\mathbb R}^d)$, 
which determines the duality $\<\cdot\,,\,\cdot\>$ between $V$ 
and $V^*$. 
Using Assumptions (A3), (A7) and (A8) with $l=r$, 
one can easily show that there exist a constant $C$ 
and a random variable $\xi$ such that $E\xi^2<\infty$ 
and 
\begin{align*}
& \Big|\sum_{|\alpha|\leq1,|\beta|\leq1} 
(-1)^{|\alpha|}\big(a^{\alpha\beta}\, D^{\beta}v\, ,\,
D^{\alpha}\varphi
\big)_r\Big|\leq C|v|_{r+1}|\varphi|_{r+1}, 
\\
& \sum_{k=1}^{d_1}
|(b^{\alpha}_kD^{\alpha}v\, ,\,
\varphi\big)_r|^2
\leq C |v|_{r+1}^2|\varphi|_r^2, 
\\
& |\big( f(t)\, ,\, \varphi\big)_{r}|^2  
\leq \xi|\varphi|^2_{r},
\quad
\sum_{k=1}^{d_1}|(g_k(t)\, ,\, \varphi\big)_r|^2
\leq  \xi|\varphi|_r^2
\end{align*}
for all  $\omega$, 
$t\in[0,T]$ and $v,\varphi\in W^r_2({\mathbb R}^d)$. 
Therefore the operators $A(t,\cdot)$, $B_k(t,\cdot)$ defined by 
\begin{align} \nonumber
\<A(t,v),\varphi\>& =
\sum_{|\alpha|\leq1,|\beta|\leq1} 
(-1)^{|\alpha|}\big(a^{\alpha\beta}\, D^{\beta}v\, ,\,
D^{\alpha}\varphi
\big)_r +\big( f(t)\, ,\, \varphi\big)_r,
\\
                                       \label{ABbis}
(B_k(t,v)\,,\,\varphi)& =
\big(b^{\alpha}_kD^{\alpha}v\, ,\,
\varphi\big)_r+ 
\big(g_k(t)\, ,\, \varphi\big)_r, 
\quad v,\varphi\in V  
\end{align}
are mappings from $V$ into $V^*$ and $H$, respectively, 
for each $k$ and $\omega,t$, such that the growth 
conditions \eqref{growthB} and \eqref{growthA} hold. 
Thus we can cast the 
Cauchy problem \eqref{11.12.06.06}--
\eqref{21.29.05.06}  into the evolution equation 
\eqref{u}, and it is 
an easy to verify that conditions (C1)--(C4) hold. 
Thus this evolution equation admits a unique 
solution $u$, which clearly a generalized solution 
to \eqref{11.12.06.06}--
\eqref{21.29.05.06}. 
Due to assumptions (A1)--(A3) and (A7)--(A8) 
by Theorem 1.1 of \cite{KrRo} 
$u$ is a $W^{r+2}({\mathbb R}^d)$-valued stochastic 
process such that 
$$
E\sup_{t\leq T}|u(t)|^2_{r+2}
+E\int_0^T|u(t)|^2_{r+3}\,dt<\infty. 
$$
Hence it is obvious that (T3) holds, 
and it is easy to verify 
(T2) with $\nu=\frac{1}{2}$ like it is done in the proof of 
 Corollary \ref{coregularity}. 
 Finally, it is an easy exercise  to show 
that (T1) holds. Now we can finish the proof 
of the theorem by applying Theorem \ref{convimplunif}.    
\end{proof}

From the previous theorem we obtain the following corollary
by Sobolev's embedding from $W^r_2$ to ${\mathcal C}^q$. 
\begin{cor}
Let $q$ be any non-negative number and assume that the assumptions of
 Theorem \ref{th4.6} hold with $r > q+\frac{d}{2}$. Then there exist
 modifications $\bar{u}$ and $\bar{u}^\tau$ of $u$ and $u^\tau$, respectively, such that
the derivatives $D^\gamma \bar{u}$ and $D^\gamma \bar{u}^\tau$  in $x$ up to order $q$ 
are functions  continuous in $x$. Moreover, there exists a constant $C$ independent of $\tau$ such that
\begin{align}                              
E\max_{0\leq i
\leq m}\sup_{x\in \RR^d}& \sum_{|\gamma|\leq q}|D^\gamma \big( \bar{u}(i\tau,x)-\bar{u}^{\tau}(i\tau,x)\big) |^2
\nonumber \\
& \quad +E\sum_{i=1}^{m} \sup_{x\in \RR^d} \sum_{|\gamma|\leq q+1} 
|D^\gamma \big( \bar{u}(i\tau,x)-\bar{u}^{\tau}(i\tau,x)\big|^2\tau
\leq C\tau .
\end{align}

\end{cor} 

\noindent {\bf Acknowledgments:} 
The authors wish to thank the referee for
helpful comments.

\end{document}